\documentclass[twoside,leqno]{article}

\usepackage{amssymb}
\usepackage{amsmath}
\usepackage[active]{srcltx}
\usepackage[cp1250]{inputenc}
\usepackage{enumerate}
\input xy
\xyoption{all}
\newcommand{\punkt}{\hspace{-0.3cm} . \hspace{0.1cm}}

\numberwithin{equation}{section}
\newtheorem{thm}[equation]{\sc Theorem}
\newtheorem{lem}[equation]{\sc Lemma}
\newtheorem{cor}[equation]{\sc Corollary}
\newtheorem{prop}[equation]{\sc Proposition}
\newtheorem{rem}[equation]{\sc Remark }

\newtheorem{defin}[equation]{\sc Definition}
\newtheorem{ex}[equation]{\sc Example}

\makeatletter
\renewcommand{\@seccntformat }[1]{\csname the#1\endcsname. }
\makeatother

 1

\def\CH{{\cal H}}

\def\CK{{\cal K}}

\def\CL{{\cal L}}

\def\CR{{\cal R}}

\def\CY{{\cal Y}}
\def\CS{{\cal S}}

\def\vsp{\vspace*{1.5ex}}

\def\epv {{$\mbox{}$\hfill ${\Box}$\vspace*{1.5ex} }}

\def\ov{\overline}

\def\bdim{\mbox{\bf dim}}
\def\ind{\mbox{{\rm ind}}}
\def\mod{\mbox{{\rm mod}}}

\def\Ext{\mbox{\rm Ext}}

\def\ya{\mbox{\em (a) \hspace{0.3em}}}
\def\yb{\mbox{\em (b) \hspace{0.3em}}}
\def\yc{\mbox{\em (c) \hspace{0.3em}}}
\def\yd{\mbox{\em (d) \hspace{0.3em}}}
\def\ye{\mbox{\em (e) \hspace{0.3em}}}
\def\yf{\mbox{\em (f) \hspace{0.3em}}}

\def\yi{\mbox{\em (i) \hspace{0.3em}}}
\def\yii{\mbox{\em (ii) \hspace{0.3em}}}
\def\yiii{\mbox{\em (iii) \hspace{0.3em}}}
\def\yiv{\mbox{\em (iv) \hspace{0.3em }}}
\def\yv{\mbox{\em (v) \hspace{0.3em}}}

\def\vsp{\vspace*{1.5ex}}

\def\vedge#1{{\buildrel{#1} \over {\hbox to
20pt{\hspace{-0.2em}$-$\hspace{-0.2em}$-$\hspace{-0.2em}$-$ }}}}

\begin{document}

\renewcommand{\thefootnote}{}
\def\refname{\begin{center}\normalsize{\it REFERENCES}\end{center}}

 \title{\large{\bf Lie algebras associated with quadratic forms and their applications to
Ringel-Hall algebras}
   \\
   \vspace{0.5em}}
\author{  \normalsize {\bf Justyna  Kosakowska }
\\
 \small{\it Faculty of Mathematics and Computer Science,
Nicolaus Copernicus University}\\ \small{\it ul. Chopina 12/18,
87-100 Toru\'n, Poland} \\
\small{\it E-mail: justus@mat.uni.torun.pl}}
\date{}
\maketitle

%\footnote{2000 {\it Mathematics Subject Classification.} Primary
%16G20, 16W10; Secondary 16G60.}

\footnote{Partially supported by Research Grant No. N N201 269 135 of Polish Ministry of Science and High Education}

\footnotesize {\bf Abstract.}  We define and investigate nilpotent Lie algebras
associated with quadratic forms. We also present their connections with
Lie algebras and Ringel-Hall algebras associated with representation directed algebras. \vsp \vsp

\noindent{\small\bf MR(2000) Subject Classification}\ \ {\rm
17B20}

\normalsize

\section{Introduction}

Let $q$ be a~unit integral quadratic  form
$$q(x)=q(x(1),\ldots,x(n))=\sum_{i=1}^nx(i)^2+\sum_{i,j}a_{ij}x(i)x(j),$$
where $a_{ij}\in \{-1,0,1\}$. In \cite{bkl}, 
with $q$ complex Lie algebras $G(q)$, $\widetilde{G}(q)$ are associated, where 
$\widetilde{G}(q)$ is the extension of $G(q)$ by the $\mathbb{C}$-dual of the
radical of $q$ and $\mathbb{C}$ is the complex number field. 

The following facts were proved in \cite{bkl}. 
\begin{itemize}
 \item If
$q$ is positive definite and connected, then
$G(q)=\widetilde{G}(q)$ is a~finite dimensional simple Lie
algebra.
\item If $q$ is connected and non-negative of corank one or
two, then $\widetilde{G}(q)$ is isomorphic to an~affine Kac-Moody algebra (if the
corank of $q$ equals one) or to elliptic (if the corank of $q$ equals
two and $q$ is not of Dynkin type $\mathbb{A}_n$).
\end{itemize}
  In \cite{bkl},
the Lie algebra $G(q)$ was defined by generators and relations.
Unfortunately, the set of relations defining $G(q)$ is infinite.
In \cite{barotrivera}, the authors give a~finite and small set of
relations sufficient to define $G(q)$ for positive definite forms
$q$.

In this paper, for any integral quadratic form $q$ (\ref{eq:form}), we define 
by generators and relations, a~Lie algebra
$L(q,\mathfrak{r})$. For a~positive definite form
$q$, we describe a~minimal set of relations defining $L(q,\mathfrak{r})$.
Moreover we show that, for any representation directed
$\mathbb{C}$-algebra $A$ with Tits form $q_A$, there are isomorphisms
of Lie algebras $$L(q_A,\mathfrak{r})\cong \CL(A) \cong \CK(A), $$
where $\CL(A)$ is the~Lie algebra associated with $A$ in \cite{riedtmann} by Ch. Riedtmann
and $\CK(A)$ is the~Lie algebra associated with $A$ in \cite{ringel1} by C. M. Ringel.
The isomorphism $\CL(A) \cong \CK(A)$ is proved in \cite{kaskos}.
Results of the present paper allow us to define Lie algebras $\CL(A)$ and
$\CK(A)$ in a~combinatorial way.  
Similar results are presented in
\cite{kos05b} and \cite{kos05c} for Tits forms of posets of finite prinjective type.

The paper is organised as follows.

  In Section \ref{sec:forms} we give basic definitions and facts concerning
       weakly positive and positive definite quadratic forms and their roots.
 
 In Sections \ref{sec:lie-alg}, \ref{sec:root-space} we give a~definition and 
prove basic properties
       of the Lie algebra $L(q,\mathfrak{r})$. 
  Moreover (for some class of quadratic form $q$) we show that the Lie algebra $L(q,\mathfrak{r})$
 is a~Lie subalgebra of $G(q)$. 
 
 In Section \ref{sec:ringel-hall} we prove the existence of isomorphisms
 of Lie algebras $$L(q_A,\mathfrak{r})\cong \CL(A) \cong \CK(A), $$ for any representation
 directed $\mathbb{C}$-algebra $A$. Moreover we present applications 
of these results to Ringel-Hall algebras of representation directed algebras.
 
 In Section \ref{sec:min-set} we give a~minimal set of relations sufficient to define
       the Lie algebra $L(q,\mathfrak{r})$, where $q$ is a~positive definite quadratic form.
   More precisely we prove the following theorem.
\begin{thm}\punkt
 Let $q$ be a~positive definite quadratic form {\rm (\ref{eq:form})}. There is an~isomorphism
of Lie algebras $$L(q,\mathfrak{r})\cong L(q)/(\mathfrak{j}), $$ where
$L(q)$ is the free complex Lie algebra generated by the set $\{v_1,\ldots,v_n\}$
and $(\mathfrak{j})$ is the ideal of $L(q)$ generated by the set $\mathfrak{j}$, which is consisted
of the following elements
\begin{itemize}
 \item $[v_i,v_j]$ for all $i,j\in\{1,\ldots,n\}$ such that
       $i<j$ and $a_{ij}\neq -1$,
 \item $[v_i,[v_i,v_j]]$ for all $i,j\in\{1,\ldots,n\}$ such that
       $i<j$ and $a_{ij}= -1$,
\item $[v_j,[v_i,v_j]]$ for all $i,j\in\{1,\ldots,n\}$ such that
       $i<j$ and $a_{ij}= -1$,
\item $[v_{i_1},\ldots,v_{i_m}]$ for all positive chordless cycles $(i_1,\ldots,i_m)$ 
$($see Section {\rm \ref{sec:min-set}} for definitions$)$.
\end{itemize}

\end{thm}

  Finally, in Section \ref{sec:lie-rep-dir}, we present some examples and remarks.

\section{Preliminaries on weakly positive and positive definite quadratic forms}\label{sec:forms}

Let $e_1,\ldots,e_n$ be the standard basis of the free abelian group $\mathbb{Z}^n$.
Let $q:\mathbb{Z}^n\to\mathbb{Z}$ be a~connected unit integral quadratic form defined
by
\begin{equation}
 q(x)=q(x(1),\ldots,x(n))=\sum_{i=1}^nx(i)^2+\sum_{i,j}a_{ij}x(i)x(j), \label{eq:form}
\end{equation}
where $a_{ij}\in \{-1,0,1\}$. Let $B(q)$ be the~bigraph associated
with $q$ (i.e. the set of vertices of $B(q)$ is $\{1,\ldots,n\}$;
for $i\neq j$ there exists a~solid edge $\xymatrix{i\ar@{-}[r] &
j}$ if and only if $a_{ij}=-1$ and a~broken edge
$\xymatrix{i\ar@{--}[r] & j}$ if and only if $a_{ij}=1$).
An~integral quadratic form $q$ is said to be \textbf{weakly
positive}, if $q(x)>0$ for any $0\neq x\in \mathbb{N}^n$, where $\mathbb{N}$ is
the set of non-negative integers.

A~vector $x\in\mathbb{Z}^n$ is called \textbf{a~root} of $q$, if
$q(x)=1$; if in addition $x(i)\geq 0$, for any $i=1,\ldots,n$,
then we call $x$ \textbf{a~positive root}. Denote by
\begin{equation}
\CR_q=\{ x\in\mathbb{Z}^n\; ;\; q(x)=1\},\;\; \CR_q^+=\{
x\in\mathbb{N}^n\; ;\; q(x)=1\} \label{eq:roots}\end{equation} the
set of all roots and all positive roots of $q$, respectively.

We associate with $q$ the~symmetric $\mathbb{Z}$-bilinear form
\begin{equation}
\langle -,-\rangle_q  :\mathbb{Z}^n\times \mathbb{Z}^n\to\mathbb{Z}, \label{eq:diff}
\end{equation} where
$\langle x,y\rangle_q=q(x+y)-q(x)-q(y)$, for all $x,y\in\mathbb{Z}$. It is
straightforward to check that
\begin{equation}
  \langle e_i,x\rangle_q=2\cdot x(i)+\sum_{i\neq j}a_{ij} x(j), \label{eq:diff-explicite}
\end{equation}

\noindent for any $i=1,\ldots,n$. Let us recall the following useful facts concerning the $\mathbb{Z}$-linear form
$\langle e_i,-\rangle_q $ (see \cite{ri}). \vsp

\begin{lem}\punkt  Let $q:\mathbb{Z}^n\to \mathbb{Z}$ be an~integral quadratic form
{\rm (\ref{eq:form})} and let $i\in \{1,\ldots,n\}$.

\ya  $\langle e_i,e_i\rangle_q=2$.

\yb Let $x$ be a~root of $q$ and $0\neq
d\in \mathbb{Z}$. The vector $x-de_i$ is a~root of $q$ if and only
if $d=\langle e_i,x\rangle_q  $.

{\rm (b$'$)} Let $x$, $y$ be roots of $q$. The vector $x+y$ is a~root of $q$
if and only if $\langle x,y\rangle_q=-1$.

Assume that $q$ is positive definite.

\yc If $x$ is a~root of $q$ such that $x\neq e_i$, then
$$-1\leq\;\langle e_i,x\rangle_q  \leq 1.$$

\yd The set $\CR_q$ of all roots of $q$ is finite.

Assume that $q$ is weakly positive.

\ye If $x$ is a~root of $q$, then
$$-1\leq\;\langle e_i,x\rangle_q  .$$

\yf The set $\CR_q^+$ of all positive roots of $q$ is finite.
\epv\label{lem:pierwiastki}\end{lem}

\begin{lem}\punkt  Let $q:\mathbb{Z}^n\to \mathbb{Z}$ be a~weakly positive 
quadratic form and let $z\neq 0$ be a~positive root of $q$. Then
$z$ is a~{\bf Weyl root}, i.e. there exists a~chain
$$x^{(1)},\ldots ,x^{(m)} $$
of positive roots of $q$ such that

\ya $x^{(1)}=z$, $x^{(i)}=x^{(i-1)}-e_{j_i}$ for $i=1,\ldots,m$
and for some $j_i\in \{1,\ldots,n\}$,

\yb $x^{(m)}=e_j$, for some $j\in \{1,\ldots,n\}$.\epv
\label{lem:weyl}\end{lem}

\section{Lie algebras associated with quadratic forms}\label{sec:lie-alg}

In a~complex Lie algebra $L$, we use the following multibracket notation
\begin{equation}[y_1,y_2,\ldots
,y_n]=[y_1,[y_2,[\ldots[y_{n-1},y_n]]]],
\label{eq:standard-form}\end{equation} for all $y_1,\ldots,y_n\in
L$. Let $L$ be a~complex Lie algebra generated by a~set 
$\{v_1,\ldots,v_n\}$. Elements of the form $[v_{i_1},\ldots,v_{i_m}]$ we call
\textbf{standard multibrackets}. All Lie algebras considered in this
paper are assumed to be complex finitely generated Lie algebras
and all quadratic forms are assumed to be integral quadratic
forms (\ref{eq:form}).

Let $S_n$ be the symmetric group of $n$-elements.

\begin{lem}\punkt  \ya Let $L$ be a~Lie algebra and let $y_1,\ldots,y_n,x\in L$.
There exists a~subset $\CS\subseteq S_n$ and for any $\sigma\in\CS$
there exists $\varepsilon_{\sigma}\in \{0,1\}$, such that
$$[[y_1,\ldots, y_n], x]=\sum_{\sigma\in\CS}(-1)^{\varepsilon_{\sigma}}
[y_{\sigma(1)},\ldots, y_{\sigma(n)},x]. $$

\yb Let $L$ be a~Lie algebra generated by a~set
$\{y_1,\ldots,y_n\}$. Any element $y\in L$ is a~linear combination
of standard multibrackets  $[y_{i_1},\ldots, y_{i_m}]$,
where $i_j\in \{1,\ldots ,n\}$. \label{lem:lie-basis}\end{lem}

{\bf Proof.} (a) Apply recursively the Jacobi identity, see also
\cite[Lemma 1.1]{as}.

The statement (b) follows from (a).\epv

With a~quadratic form $q$ (\ref{eq:form}), we associate the
complex free Lie algebra
\begin{equation}
 L(q)={\rm Lie}_{\mathbb{C}}\langle v_1,\ldots ,v_n\rangle \label{eq:lie}
\end{equation}
generated by the set $\{v_1,\ldots,v_n\}$.
Note that the Lie algebra $L(q)$ has a~$\mathbb{N}^n$-gradation,
if we define the degree of $v_i$ to be $e_i$, for any
$i=1,\ldots,n$.

Let $\mathfrak{a}\subseteq L(q)$ be a~subset consisting of some
standard multibrackets and let $(\mathfrak{a})$ be
the~homogeneous ideal of the Lie algebra $L(q)$ generated by the
set $\mathfrak{a}$. Let
\begin{equation} L(q,\mathfrak{a})=L(q)/(\mathfrak{a})\label{eq:lie-ideal}\end{equation} be the quotient
Lie algebra with induced $\mathbb{N}^n$-grading. Denote by $\pi:L(q)\to
L(q,\mathfrak{a})$ the natural epimorphisms. Let
$v=[v_{i_1},\ldots,v_{i_m}]\in L(q)$. For the sake of simplicity,
we will denote by
$v=[v_{i_1},\ldots,v_{i_m}]$ the element $\pi(v)$.

$\bullet$ We call an~element $v=[v_{i_1},\ldots ,v_{i_m}]\in
L(q)$ a~{\bf root}, if $m\geq 1$ and $\langle
e_{i_k},e_{i_{k+1}}+\ldots +e_{i_m}\rangle_v  =-1$, for all
$k=1,\ldots ,m-1$.

%$\bullet$ Let $v=[v_{i_1},\ldots ,v_{i_m}]\in L(q)$.
%We call the set ${\rm supp}(v)=\{i_1,\ldots ,i_m\}\subseteq
%\{1,\ldots,n\}$ the~{\bf support} of the element $v$.

$\bullet$ Let $v=[v_{i_1},\ldots ,v_{i_m}]\in L(q)$, we set
$\ell(v)=m$ and we call this number the~\textbf{length} of the element $v$.

$\bullet$ We set $e_v=e_{i_1}+\ldots +e_{i_m}\in\mathbb{N}^n$.

For $e\in\mathbb{N}^n$,
denote by $L(q,\mathfrak{a})_e$ the homogeneous space spanned by all
standard multibrackets  $v=[v_{i_1},\ldots ,v_{i_m}]\in
L(q,\mathfrak{a})$, of degree $e_v=e$. Moreover, for any integer $m$, let
\begin{equation}
 L(q,\mathfrak{a})_m=\displaystyle\bigoplus_{\substack{e\in\mathbb{N}^n \\
 e(1)+\ldots +e(n)\leq m}}L(q,\mathfrak{a})_e \;\;\mbox{ and }\;\; (\mathfrak{a})_m
= (\mathfrak{a})\cap L(q)_m.
\end{equation}

The following lemma shows connections between roots in the complex Lie
algebra $L(q)$ and positive roots of the quadratic form $q$.

\begin{lem}\punkt Let $u,w,v=[v_{i_1},\ldots ,v_{i_m}]\in L(q)$ be roots.

\ya At least one of the elements $[v_{i_1},v_{i_2}]$,
$[v_{i_2},v_{i_1},v_{i_3},\ldots ,v_{i_m}]$ is not a~root.

\yb If $\langle e_u,e_v+e_w\rangle_q=-1$,
then $\langle e_u,e_v\rangle_q\geq 0$ or $\langle e_u,e_w\rangle_q\geq 0$.

\yc The vector $e_v=e_{i_1}+\ldots +e_{i_m}$ is a~positive root of
the quadratic form $q$.

Let $v=[v_{i_1},\ldots ,v_{i_m}]\in L(q)$ and let $q$ be weakly positive.

\yd If $e_v=e_{i_1}+\ldots +e_{i_m}$ is a~positive root of $q$,
then there exists a~permutation $\sigma\in S_m$ such that the
element $v_{\sigma}=[v_{i_{\sigma(1)}},\ldots ,v_{i_{\sigma(m)}}]$
is a~root in $L(q)$. \label{lem:root} \end{lem}

{\bf Proof.} (a) Let $v=[v_{i_1},\ldots ,v_{i_m}]\in L(q)$ be
a~root and  assume that $[v_{i_1},v_{i_2}]$ is a~root. It follows
that $\langle e_{i_1},e_{i_2}\rangle_q  =-1$ and $\langle
e_{i_1},e_{i_3}+\ldots + e_{i_m}\rangle_q  =\langle
e_{i_1},e_{i_2}+\ldots +e_{i_m}\rangle_q  -\langle
e_{i_1},e_{i_2}\rangle_q  =-1-(-1)=0$. Therefore
$[v_{i_2},v_{i_1},v_{i_3},\ldots ,v_{i_m}]$ is not a~root.

The statement (b) is obvious.

(c) Obviously $e_{i_m}$ is a~root of $q$. Let $2\leq k\leq m$ and
assume that $e_{i_k}+\ldots +e_{i_m}$ is a~root of $q$. We have
$\langle e_{i_{k-1}},e_{i_k}+\ldots +e_{i_m}\rangle_q=-1$, because
$v$ is a~root in $L(q)$.  By Lemma
\ref{lem:pierwiastki}(b), the vector
$e_{i_{k-1}}+e_{i_k}+\ldots+e_{i_m}$ is a~root of $q$. Inductively
we finish the proof of the statement (c).

(d) Let $v$ be a~positive root of $q$. From Lemma
\ref{lem:pierwiastki}(b) it follows that  $v+e_i$ is a~positive
root of $q$ if and only if $\langle e_i,v\rangle_q  =-1$.
Therefore the statement (d) follows easily from Lemma
\ref{lem:weyl} and Lemma \ref{lem:pierwiastki}(b), because $q$ is
weakly positive.\epv

\begin{defin}\punkt {\rm (a) Let $\mathfrak{r}$ be the set of all standard multibrackets
$[v_{i_1},\ldots ,v_{i_m}]\in L(q)$, such that $[v_{i_1},\ldots
,v_{i_m}]$ is not a~root and $[v_{i_2},\ldots ,v_{i_m}]$ is
a~root, where $i_1,\ldots,i_m\in \{1,\ldots,n\}$ and
$m\in\mathbb{N}$.

(b) Let $$L(q,\mathfrak{r})=Lie_{\mathbb{C}}\langle
v_1,\ldots,v_n\rangle /(\mathfrak{r})
$$ be a~complex Lie algebra generated by the set $\{v_1,\ldots,v_n\}$
modulo the ideal $(\mathfrak{r})$ generated by the set
$\mathfrak{r}$. We consider $L(q,\mathfrak{r})$ as a~Lie algebra with
a~$\mathbb{N}^n$-gradation, where we define the degree of $v_i$ to
be $e_i$, for any $i=1,\ldots,n$.} \label{def:lie}\end{defin}

\begin{lem}\punkt  Let $\mathfrak{a}\subseteq L(q)$ be a~subset consisting of some standard 
multibrackets.
Let $v=[v_{i_1},\ldots,v_{i_m}]$, $y\in L(q,\mathfrak{a})$.
Assume that for any $z\in L(q,\mathfrak{a})$, such that
$\ell(z)\leq \ell(v)+\ell(y)$, the following condition is satisfied:

\begin{equation} \mbox{ if } \langle e_i,e_z\rangle_q\neq -1,
\mbox{ then } [v_i,z]=0 \mbox{
in } L(q,\mathfrak{a}).\label{eq:zero}\end{equation}

\noindent Then $[v,y]=0$ in $L(q,\mathfrak{a})$ or there
exists $\sigma \in S_m$ and $\varepsilon\in \{0,1\}$ such that
$$[v,y]=(-1)^\varepsilon[v_{i_{\sigma(1)}},\ldots,v_{i_{\sigma(m)}},y]$$
in $L(q,\mathfrak{a})$.
 \label{lem:root1}\end{lem}

{\bf Proof.} Let $v=[v_{i_1},\ldots,v_{i_m}],y\in
L(q,\mathfrak{a})$. We precede with induction on $m=\ell(v)$. For
$m=1$, the lemma is obvious.

Let $m>1$. We apply the Jacobi identity and get
$$[v,y]=-[y,v]=[v_{i_1},[[v_{i_2},\ldots,v_{i_m}],y]]-[[v_{i_2},\ldots,v_{i_m}],[v_{i_1},y]]. $$
Note that $$\langle
e_{i_1},e_{i_2}+\ldots+e_{i_m}+e_y\rangle_q=\langle
e_{i_1},e_{i_2}+\ldots+e_{i_m}\rangle_q+\langle
e_{i_1},e_y\rangle_q.$$
Therefore
$\langle e_{i_1},e_{i_2}+\ldots+e_{i_m}\rangle_q\neq -1$ or $\langle e_{i_1},e_y \rangle_q\neq -1$ or $\langle
e_{i_1},e_{i_2}+\ldots+e_{i_m}+e_y\rangle_q\neq -1$.

If $\langle e_{i_1},e_{i_2}+\ldots+e_{i_m}\rangle_q\neq -1$, then, by (\ref{eq:zero}),
$v=[v_{i_1},\ldots,v_{i_m}]=0$ and $[v,y]=0$. We are done.
If $\langle
e_{i_1},e_y\rangle_q\neq -1$, then, by (\ref{eq:zero}), we have $[v_{i_1},y]=0$ and
$[v,y]=[v_{i_1},[[v_{i_2},\ldots,v_{i_m}],y]]$. We finish by
induction on $m$ applied to $[[v_{i_2},\ldots,v_{i_m}],y]]$.
In the case $\langle e_{i_1},e_{i_2}+\ldots+e_{i_m}+e_y\rangle_q\neq
-1$, we have
$[v_{i_1},[[v_{i_2},\ldots,v_{i_m}],y]]=0$.
Finally
$$[v,y]=(-1)^{\varepsilon}[[v_{i_2},\ldots,v_{i_m}],[v_{i_1},y]]$$ and we finish by induction on $m$. \epv

\section{Grading of $L(q,\mathfrak{r})$ and positive roots of $q$}\label{sec:root-space}

\begin{lem}\punkt Let $q$ be a~weakly positive quadratic form {\rm (\ref{eq:form})} and
let $\mathfrak{a}\subseteq L(q)$ be a~subset consisting of some standard
multibrackets. Let $m$ be a~positive integer. Assume that
the following conditions are satisfied.

\ya If $v=[v_{i_1},\ldots v_{i_s}]$ is not a~root in $L(q)$ and $\ell(v)=s\leq m$, then
$$L(q,\mathfrak{a})_{e_v}=0.$$

\yb If $v=[v_{i_1},\ldots v_{i_s}]$ is a~root in $L(q)$ and $\ell(v)=s< m$, then $$\dim_{\mathbb{C}}L(q,\mathfrak{a})_{e_v}\leq 1.$$

\noindent Then $\dim_{\mathbb{C}}L(q,\mathfrak{a})_{e_v}\leq 1$, if $v=[v_{i_1},\ldots v_{i_s}]$
is a~root in $L(q)$ and $\ell(v)=s= m$.

\label{lem:root-spaces}\end{lem}

\textbf{Proof.} Let $v=[v_{i_1},\ldots ,v_{i_m}]$ be a~root of
$L(q)$. If $v=0$ in $L(q,\mathfrak{a})$, then we
are done. Assume that $0\neq v\in
L(q,\mathfrak{a})$. We have to prove that
$\dim_{\mathbb{C}}L(q,\mathfrak{a})_{e_v}\leq 1$. If $m=1$, the
lemma is obvious. Let $m>1$, and let $0\neq w\in
L(q,\mathfrak{a})_{e_v}$ be a~standard multibracket. It follows that there
exists a~permutation $\sigma\in S_m$ such that
$w=v_\sigma=[v_{i_{\sigma(1)}},\ldots ,v_{i_{\sigma(m)}}]\neq 0$.
Since $v_{\sigma}\neq 0$ in $L(q,\mathfrak{a})$,
$\ell(v_{\sigma})=m$, then the assumption (a) of our lemma yields
that $v_{\sigma}$ is a~root. It is enough to prove that there
exists $a\in \mathbb{C}$ such that $v_\sigma=av$. Since
$v,v_{\sigma}\in L(q,\mathfrak{a})_{e_v} $, there exists
$k=1,\ldots ,m$ such that $i_k=i_{\sigma(1)}$.  Note that we may
assume, without loss of the generality, that $k<m$, because
$[v_{i_{m-1}},v_{i_m}]=-[v_{i_m},v_{i_{m-1}}]$ and we may replace
$v$ by $-v$.

If $k=1$, then
$\overline{v}=[v_{i_{\sigma(2)}},\ldots ,v_{i_{\sigma(m)}}]$, $[v_{i_2},\ldots
,v_{i_m}]\in L(q,\mathfrak{a})_{e_{\overline{v}}}$, $\ell({\overline{v}})<m$, and 
the condition (b) yields
that $\dim_{\mathbb{C}}L(q,\mathfrak{a})_{e_{\overline{v}}}\leq 1$. Then
there exists $a\in \mathbb{C}$ such that
$[v_{i_{\sigma(2)}},\ldots ,v_{i_{\sigma(m)}}]=a[v_{i_2},\ldots
,v_{i_m}]$. Therefore $v_\sigma=av$ and we are done.

Let $k>1$. Consider the following set
$$\CY=\{v_{\tau}\; ; \; \mbox{ for all } \tau\in S_m
\mbox{ such that there exists } c\in \mathbb{C} \mbox{ such that }
v_\tau=cv \}.
$$ Note that for all $v_{\tau}\in \CY$ there exists $l$ such that
$i_{\tau(l)}=i_{\sigma(1)}$. We choose an element $v_\tau\in \CY$
such that $l$ is minimal with this property. Without loss
of the generality, we may assume that $\tau={\rm id}$, $v_\tau=v$ and $k=l$.
Let $$v=[v_{i_1},[v_{i_2},[\ldots,v_{i_{k-1}},[v_{i_k},y]\ldots]]]
=[v_{i_1},v_{i_2},\ldots,v_{i_{k-1}},v_{i_k},y],$$
where we set $y=[v_{i_{k+1}},\ldots,v_{i_m}]$.
By the choice of $k$ it follows that $i_j\neq
i_{\sigma(1)}$ for all $j< k$. Our assumptions yield: $\langle
e_{i_k},e_y\rangle_q  =-1$, because $v$ is a~root and $$\langle
e_{i_k},e_y+e_{i_1}+\ldots
+e_{i_{k-1}}\rangle_q  =\langle
e_{i_{\sigma(1)}},e_{i_{\sigma(2)}}+\ldots
+e_{i_{\sigma(m)}}\rangle_q  =-1,$$ because $v_{\sigma}$ is a~root. 
By the bilinearity of $\langle -,-\rangle_q$, it follows that $\langle
e_{i_k},e_{i_1}+\ldots +e_{i_{k-1}}\rangle_q
=0$. We prove that $$v=[[v_{i_1},\ldots
,v_{i_k}],y]$$ in $L(q,\mathfrak{a})$.

Applying the Jacobi identity, we get
\footnotesize $$\begin{array}{l}[v_{i_1},\ldots ,v_{i_{k-1}},[v_{i_k},
y]]= -[v_{i_1},\ldots
,v_{i_{k-2}},[v_{i_k},[y, v_{i_{k-1}}]]]-
[v_{i_1},\ldots ,v_{i_{k-2}},[y,[v_{i_{k-1}},
v_{i_k}]]].\end{array}$$ \normalsize Note that
$$\langle e_{i_{k-1}},e_y\rangle_q  =\langle e_{i_{k-1}},e_y+e_{i_k}\rangle_q
 -\langle e_{i_{k-1}},e_{i_k}\rangle_q
= -1-\langle e_{i_{k-1}},e_{i_k}\rangle_q  .$$ If
$\langle e_{i_{k-1}},e_y\rangle_q  =-1$, then $\langle
e_{i_{k-1}},e_{i_k}\rangle_q  =0$. Therefore the condition (a) yields $[v_{i_{k-1}},v_{i_k}]=0$ and
$$v=[v_{i_1},\ldots
,v_{i_{k-2}},[v_{i_k},[v_{i_{k-1}}, y]]]$$
which is the contradiction with the choice of $k$.
Therefore
$$v=[v_{i_1},\ldots ,v_{i_{k-2}},[[v_{i_{k-1}},
v_{i_k}]\, y]].$$ Inductively, applying the Jacobi identity to
$$[v_{i_1},\ldots
,v_{i_s},[[v_{i_{s+1}},\ldots, v_{i_k}], y]]
$$ we get \small
$$\begin{array}{l} [v_{i_1},\ldots ,v_{i_s},[[v_{i_{s+1}},\ldots, v_{i_k}],
y]]=\\ -[v_{i_1},\ldots
,v_{i_{s-1}},[[v_{i_{s+1}},\ldots,v_{i_k}],[y,v_{i_s}]]]
-  [v_{i_1},\ldots
,v_{i_{s-1}},[y,[v_{i_s},[v_{i_{s+1}},\ldots,v_{i_k}]]]].
\end{array}
$$ \normalsize Consider
$$\begin{array}{c}
\langle e_{i_s},e_y\rangle_q  =\langle
e_{i_s},e_y+e_{i_{s+1}}+\ldots+e_{i_k}\rangle_q
-\langle
e_{i_s},e_{i_{s+1}}+\ldots+e_{i_k}\rangle_q
=
\\ = -1-\langle e_{i_s},e_{i_{s+1}}+\ldots+e_{i_k}\rangle_q  .\end{array}$$ If
$\langle e_{i_s},e_y\rangle_q  =-1$, then $\langle
e_{i_s},e_{i_{s+1}}+\ldots+e_{i_k}\rangle_q
=0$. The condition (a) yields $[v_{i_s},[v_{i_{s+1}},\ldots,v_{i_k}]]=0$ and
$$v=[v_{i_1},\ldots,
v_{i_{s-1}},[[v_{i_{s+1}},\ldots,v_{i_k}],[v_{i_s},
y]]].$$  Applying Lemma \ref{lem:root1}, we get the~contradiction with the
choice of $k$. Therefore
$$v=[v_{i_1},\ldots
,v_{i_{s-1}},[[v_{i_s},[v_{i_{s+1}},\ldots,v_{i_k}]],y]].$$
Inductively we get $v=[[v_{i_1},\ldots,v_{i_k}],y]$.

Since $v\neq 0$, $k<m$, then (by the assumption (a) of the lemma) the element  $[v_{i_1},\ldots
,v_{i_k}]$ is a~root. Therefore, by Lemma \ref{lem:root} (c), we have
$q(e_{i_1}+\ldots +e_{i_k})=1$. Now consider
$$1=q(e_{i_1}+\ldots
+e_{i_k})=q(e_{i_1}+\ldots
+e_{i_{k-1}})+q(e_{i_k})+\langle
e_{i_k},e_{i_1}+\ldots
+e_{i_{k-1}}\rangle_q.$$ Since we proved above that $\langle
e_{i_k},e_{i_1}+\ldots +e_{i_{k-1}}\rangle_q
=0$, we have
$$q(e_{i_1}+\ldots +e_{i_{k-1}})=1-q(e_{i_k})=1-1=0.$$ We get a~contradiction, because
$q$ is weakly positive. This shows that $k=1$ and $v_\sigma\in
{\cal Y}$. This finishes the proof of lemma. \epv

Let $L(q,\mathfrak{r})$ be the Lie algebra introduced in Definition
\ref{def:lie}

\begin{prop}\punkt Let $q$ be a~weakly positive quadratic form
{\rm (\ref{eq:form})}. The following conditions hold.

\ya If $e=e_{i_1}+\ldots +e_{i_m}$ is not a~root of $q$,
then $L(q,\mathfrak{r})_e=0$.

\yb If $e=e_{i_1}+\ldots +e_{i_m}$ is a~root of $q$,
then $\dim_{\mathbb{C}}L(q,\mathfrak{r})_e\leq 1$.\label{prop:root-spaces}\end{prop}

\textbf{Proof.}
(a) Assume that $e=e_{i_1}+\ldots +e_{i_m}$ is not a~root of $q$. 
By Lemma \ref{lem:root},
$v=[v_{i_1},\ldots ,v_{i_m}]$ is not a~root in $L(q,\mathfrak{r})$. Let $k$ be maximal with the property
that $[v_{i_k},\ldots ,v_{i_m}]$ is a~root. Since $v=[v_{i_1},\ldots ,v_{i_m}]$
is not a~root, then $k>1$. Therefore $[v_{i_{k-1}},v_{i_k},\ldots ,v_{i_m}]$
is not a~root and $[v_{i_{k-1}},v_{i_k},\ldots ,v_{i_n}]\in\mathfrak{r}$. Finally
$[v_{i_{k-1}},v_{i_k},\ldots ,v_{i_m}]=0$ and $v=0$ in $L(q,\mathfrak{r})$.

The statement (b) follows easily by induction on the length $\ell(v)$ of $v$, 
if we apply Lemma \ref{lem:root-spaces}.
\epv

As a~consequence we get the following corollary.

\begin{cor}\punkt Let $q$ be a~weakly positive quadratic form and let $\CR_q^+$
be the set of all positive roots of $q$. Then
$L(q,\mathfrak{r})$ is a~nilpotent Lie algebra and $$L(q,\mathfrak{r})=\bigoplus_{e\in \CR_q^+}L(q,\mathfrak{r})_e
\;\;\mbox{ and }\;\; \dim_{\mathbb{C}}L(q,\mathfrak{r})\leq
|\CR_q^+|.$$ \epv \label{cor:root-space}\end{cor}

Let $q:\mathbb{Z}^n\to \mathbb{Z}$ quadratic form (\ref{eq:form}).
With $q$ we associate Cartan matrix
$C=(c_{ij})\in\mathbb{M}_n(\mathbb{Z})$ defined by
$c_{ij}=q(e_i+e_j)-q(e_i)-q(e_j)$. Following \cite{bkl}, to $q$ we attach
a~$\mathbb{Z}^n$-graded complex Lie algebra $G(q)$ with generators
$x_i,x_{-i},h_i$, $i=1,\ldots,n$, which are homogeneous of degree
$e_i,-e_i,0$, respectively, and subject to the following
relations:
\begin{enumerate}
\item $[h_i,h_j]=0$, for all $i,j=1,\ldots ,n$,
\item $[h_i,x_{\varepsilon j}]=\varepsilon c_{ij}x_{\varepsilon
j}$, for all $i,j=1,\ldots ,n$ and $\varepsilon\in \{-1,1\}$,
\item $[x_{\varepsilon i},x_{-\varepsilon i}]=\varepsilon h_i$, for
all $i=1,\ldots ,n$ and $\varepsilon \in \{-1,1\}$,
\item $[x_{\varepsilon_1i_1},\ldots ,x_{\varepsilon_ni_n}]=0$, if
$q(\varepsilon_1e_{i_1}+\ldots +\varepsilon_ne_{i_n})>  1$ for
$\varepsilon_j\in \{-1,1\}$.
\end{enumerate}

Denote by $G^+(q)$ a~Lie subalgebra of $G(q)$
generated by the elements $x_1,\ldots,x_n$.

\begin{prop}\punkt  If $q$ is weakly positive and positive semi-definite, then
$$L(q,\mathfrak{r})\simeq G^+(q).$$ \label{prop:nonneg}\end{prop}

{\bf Proof.} By \cite[Proposition 2.2]{bkl} and Corollary \ref{cor:root-space}, we have
$$\dim_{\mathbb{C}}G^+(q)\geq
|\CR_q^+|\geq \dim_{\mathbb{C}}L(q,\mathfrak{r}).$$ On the other hand, it
is easy to see that all relations $\mathfrak{r}$ are satisfied in
$G^+(q)$. Therefore we may define a~homomorphism of Lie algebras
$$\Psi:L(q,\mathfrak{r})\to G^+(q)$$ by $\Psi(u_i)=x_i$ for all $i=1,\ldots,n$. 
Since $G^+(q)$ is generated by the elements $x_1,\ldots,x_n$,
the homomorphism $\Psi$ is surjective. Therefore $\Psi$ is
an~isomorphism, because $\dim_{\mathbb{C}}G^+(q)\geq
\dim_{\mathbb{C}}L(q,\mathfrak{r})$. \epv

\section{Connections with Ringel-Hall algebras}\label{sec:ringel-hall}

We present applications of Lie algebras $L(q,\mathfrak{r})$ to Lie algebras
and Ringel-Hall algebras associated with representation directed algebras.
We get a~description of these Ringel-Hall algebras by generators and
relations.
For the basic concepts of representation theory the reader is referred to \cite{ASS},
\cite{ARS} and for the basic concepts of Ringel-Hall algebras to \cite{ringel1},
\cite{ringel89}. 

Let $Q=(Q_0,Q_1)$ be a~finite quiver without oriented cycles. Let
$\mathbb{C}Q$ be the complex path algebra of $Q$. Assume that $I$
is an admissible ideal of $\mathbb{C}Q$ such that
$A=\mathbb{C}Q/I$ is a~representation directed algebra. 
%Let $\CR$
%be a~minimal set of elements generating the ideal $I$. 
%Thanks to
%the multiplicative basis theorem we may assume that $\CR$ is consisted
%of zero-relations and commutativity relations (see \cite{Bo1}). 
By $\mod(A)$ we denote the category of all right finite dimensional
$A$-modules and by $\ind(A)$ we denote the set of all
representatives of isomorphism classes of indecomposable
$A$-modules. For any $A$-module $M$ denote by $\bdim\,M\in\mathbb{N}^{Q_0}$ the
dimension vector of $M$ (i.e. $(\bdim\,M)(i)$ equals the number
of composition factors of $M$ which are isomorphic to the simple 
$A$-module $S_i$ corresponding to the vertex $i\in Q_0$). 
Let $q_A:\mathbb{Z}^{Q_0}\to \mathbb{Z}$ be the Tits
form of $A$ (see \cite{bo83}). By \cite[Theorem 3.3]{bo83}, $q_A$ is weakly
positive and there is a~bijection
(given by $\bdim$) between the set $\ind (A)$ and the set $\CR_{q_A}^+$. Let $\CK(A)$ be the corresponding complex Lie
algebra defined in \cite{ringel1}. Recall that, for
a~representation directed algebra $A$, the $\mathbb{C}$-Lie
algebra $\CK(A)$ is the free $\mathbb{C}$-linear space with basis
$\{u_X\; ;\; X\in\ind(A)\}$. If $X$, $Y$ are non-isomorphic
indecomposable $A$-modules such that $\Ext_A^1(X,Y)=0$, then the
Lie bracket in $\CK(A)$ is defined by
$$
[u_Y,u_X]=\left\{\begin{array}{ll} \varphi_{YX}^Z(1)\cdot u_Z &
\mbox{if there is an~indecomposable } A-\mbox{module } Z\\
&\mbox{and a~short exact sequence } \\ & 0\to X\to Z\to Y\to 0,
\\ & \\ 0 & \mbox{otherwise,}
\end{array} , \right.
 $$  where $\varphi_{YX}^Z$
are Hall polynomials (see \cite{ringel1}). In \cite{kaskos} it is
proved that the Lie algebra $\CK(A)$ is isomorphic to $\CL(A)$,
where $\CL(A)$ is the Lie algebra associated with $A$ in
\cite{riedtmann}. Let $\CH(A)$ be the universal enveloping
algebra of the Lie algebra $\CK(A)$. Recall that $\CH(A)=\CH_1(A)$, where
$\CH_q(A)$ is the generic Ringel-Hall algebra associated in \cite{ringel1} with the 
algebra $A$. In fact, in \cite{ringel1}, generic Ringel-Hall algebras were
associated with directed Auslander-Reiten quivers. However, it is possible
to associate generic Ringel-Hall algebras with representation directed
$\mathbb{C}$-algebras. The reader is referred to
\cite{kaskos} for details.

\begin{prop}\punkt Let $A$ be a~representation directed
$\mathbb{C}$-algebra.

\ya The Lie algebra $\CK(A)$ is generated by the set $\{u_i\; ;\;
i\in Q_0 \}$, where $u_i=u_{S_i}$ and $S_i$ is a~simple $A$-module
corresponding to the vertex $i\in Q_0$.

\yb In the Lie algebra $\CK(A)$ the relations from the set
$\mathfrak{r}$ hold, if we interchange $u_i$'s by $v_i$'s.
\label{prop:lie-ringel}\end{prop}

\textbf{Proof.} The statement (a) is proved in
\cite[Proposition 6]{ringel89}. Let $[v_{i_1},\ldots, v_{i_n}]$ be
an~element from the set $\mathfrak{r}$. It follows that $[v_{i_2},\ldots, v_{i_n}]$
is a~root and the element $[v_{i_1},v_{i_2},\ldots, v_{i_n}]$ is not a~root. By Lemma \ref{lem:root},
the vector $m=e_{i_2}+\ldots+e_{i_m}$ is a~positive root of the Tits form
$q_A$ of $A$.
If $[v_{i_2},\ldots, v_{i_n}]=0$ in $\CK(A)$, then we are done. Otherwise
$[v_{i_2},\ldots, v_{i_n}]=a\cdot u_M$ for some $0\neq a\in\mathbb{C}$ and
the unique indecomposable $A$-module $M\in\ind(A)$ with $\bdim\, M=m=e_{i_2}+\ldots+e_{i_n}$, because $A$ is representation directed
$\mathbb{C}$-algebra. Since $q_A$ is weakly positive, then, by Lemma \ref{lem:root}, the  vector
$e_{i_1}+m$ is not a~root of $q_A$. Therefore there exists no indecomposable $A$-module
with dimension vector $e_{i_1}+m$.
Then, by
\cite[Theorem 2]{ringel1}, $[v_{i_1},\ldots, v_{i_n}]=0$ in $\CK(A)$ and
we are done.  \epv

\begin{cor}\punkt If $A$ is a~representation directed $\mathbb{C}$-algebra,
then there is an~isomorphism of $\mathbb{C}$-algebras
  $$ F:L(q_A,\mathfrak{r})\to \CK(A)\cong \CL(A) $$
  given by $F(v_i)=u_i$, in particular $\dim_{\mathbb{C}}L(q_A,\mathfrak{r})=
|\CR_q^+|$.

If, in addition, $q_A$ is positive semi-definite, then $L(q_A,\mathfrak{r})\cong G^+(q_A)$.
\label{cor:iso-lie}\end{cor}

\textbf{Proof.} By Proposition \ref{prop:lie-ringel}, $F$ is
a~well-defined homomorphism of graded Lie algebras. 
Since the Lie algebra
$\CL(q_A,\mathfrak{r})$ is generated by the set $\{v_i\; ;\; i\in
Q_0 \}$ and the Lie algebra $\CK(A)$ is generated by the set
$\{u_i\; ;\; i\in Q_0 \}$, the homomorphism $F$ is an~epimorphism. 
By Corollary
\ref{cor:root-space}, $F$ is a~monomorphism, because $\dim_{\mathbb{C}}\CK(A)=|\CR_{q_A}^+|$.
Finally
$F$ is an~isomorphism of Lie algebras. 

The final assertion follows from Proposition \ref{prop:nonneg}.\epv

\section{A~minimal set of relations defining $L(q,\mathfrak{r})$ for a~positive definite form $q$}\label{sec:min-set}

The set $\mathfrak{r}$ usually is not a~minimal set generating the ideal $(\mathfrak{r})$ of
the Lie algebra $L(q)$. In this section we describe a~minimal set of elements defining the ideal $(\mathfrak{r})$ of $L(q)$ for a~positive definite form $q$ (\ref{eq:form}). 
In this section all quadratic forms are assumed to be positive definite.

\begin{rem}\punkt  {\em  The following easily verified facts are essentially used in this section.
\begin{enumerate}
 \item Let $i,j\in\{1,\ldots,n\}$ be such that $\langle e_i,e_j\rangle_q\neq -1$, then
$$[\ldots ,[v_i,[v_j,\ldots]]]=[\ldots ,[v_j,[v_i,\ldots]]]$$
in $L(q,\mathfrak{r})$.
Indeed, apply the Jacobi identity and note that in this case
$[v_i,v_j]\in(\mathfrak{r})$.

\item If $a\in L(q)$ is a~standard multibracket such that $e_a$ is not a~root of $q$, then $a\in (\mathfrak{r})$.
Indeed, apply Lemma \ref{lem:lie-basis} (b) and Proposition \ref{prop:root-spaces} (a).

\item Let $a,b\in L(q)$ be standard multibrackets such that $\langle e_a,e_b\rangle_q\geq 0$, then
$[a,b]\in(\mathfrak{r})$. Indeed,  $$q(e_a+e_b)=q(e_a)+q(e_b)+\langle e_a,e_b\rangle_q\geq 1+1=2, $$
then $e_a+e_b$ is not a~root of $q$.
Therefore $[a,b]\in (\mathfrak{r})$.

\item Let $a,b\in L(q)$ be standard multibrackets such that $\langle e_a,e_b\rangle_q\geq 0$, then
$$[\ldots ,[a,[b,\ldots]]]=[\ldots ,[b,[a,\ldots]]]$$
in $L(q,\mathfrak{r})$. Indeed, apply the Jacobi identity and the fact that in this case
$[a,b]\in(\mathfrak{r})$.

\item If $a\in L(q)$ and $\langle e_i,e_a\rangle_q\leq -2$, for some $i=1,\ldots,n$, 
then $a\in(\mathfrak{r})$.
Indeed, by Lemma \ref{lem:pierwiastki}, $e_a$ is not a~root of $q$ and therefore $a\in(\mathfrak{r})$.

\end{enumerate}

 }\label{remark}
\end{rem}

\subsection{The first step of reduction}

Let $\mathfrak{r}_1\subseteq\mathfrak{r}\subseteq L(q)$ be the set consisting of the following elements:
\begin{itemize}
 \item $[v_i,v_j]$ for all $i,j\in\{1,\ldots,n\}$ such that
       $i<j$ and $\langle e_i,e_j \rangle_q\neq -1$,
 \item $[v_i,[v_i,v_j]]$ for all $i,j\in\{1,\ldots,n\}$ such that
       $i<j$ and $\langle e_i,e_j \rangle_q= -1$,
\item $[v_j,[v_i,v_j]]$ for all $i,j\in\{1,\ldots,n\}$ such that
       $i<j$ and $\langle e_i,e_j \rangle_q= -1$.
\end{itemize}

\noindent Let $\mathfrak{r}_0\subseteq\mathfrak{r}\subseteq L(q)$ be the set consisting of
all elements $[v_{i_1},\ldots,v_{i_m}]$ of $\mathfrak{r}$ such that $\langle e_{i_1},e_{i_2}+\ldots +e_{i_m}\rangle_q =0$.
Define $\mathfrak{p}\subseteq \mathfrak{r}$ to be
 \begin{equation}\mathfrak{p}=\mathfrak{r}_1\cup \mathfrak{r}_0. \label{eq:p}\end{equation}

\begin{prop}\punkt If $q$ is a~positive definite quadratic form {\rm (\ref{eq:form})},
then the ideals $(\mathfrak{p})$ and $(\mathfrak{r})$
of the Lie algebra $L(q)$ are equal.  \label{prop:first-red}\end{prop}

\textbf{Proof.} The inclusion $(\mathfrak{p})\subseteq
(\mathfrak{r})$ is obvious. It is enough to prove that
$(\mathfrak{r})\subseteq (\mathfrak{p})$. Let
$v=[v_{i_2},\ldots,v_{i_m}]$ be a~root in $L(q)$ and let $i_1\in
\{1,\ldots,n\}$ be such that $[v_{i_1},v_{i_2},\ldots,v_{i_m}]\in
\mathfrak{r}$ and $\langle e_{i_1},e_{i_2}+\ldots+e_{i_m}
\rangle_q\neq 0$. Since $[v_{i_1},v_{i_2},\ldots,v_{i_m}]$ is not
a~root, we have $\langle e_{i_1},e_{i_2}+\ldots+e_{i_m}
\rangle_q\geq 1$. We claim that $[v_{i_1},v]\in \mathfrak{p}$.

We precede with induction on $\ell(v)=m-1$. For $\ell(v)< 3$ our
statement easily follows by a~case by case inspection on all
possible cases.

Let $\ell(v)\geq 3$, then $m\geq 4$, $v=[v_{i_2},[v_{i_3},
\overline{v}]]$, $\langle
e_{i_2},e_{i_3}+e_{\overline{v}}\rangle_q  =-1$, $\langle
e_{i_3},e_{\overline{v}}\rangle_q  =-1$ and
$\ell(\overline{v})\geq 1$. Note that $\langle e_{i_1},e_{i_2}\rangle_q=a_{ij}\in\{-1,0,1\}$,
if $i_1\neq i_2$. Therefore it is
enough to consider the following three cases. \vsp

{\bf 1)} If $\langle e_{i_1},e_{i_2} \rangle_q\in\{0,1\}$, then
$\langle e_{i_1},e_{i_3}+\ldots +e_{i_m}\rangle_q\geq 0$. Therefore
$[v_{i_1},v_{i_2}]$, $[v_{i_1},[v_{i_3}\ldots,v_{i_m}]]\in \mathfrak{r}$
and by the induction hypothesis we have
$[v_{i_1},[v_{i_3}\ldots,v_{i_m}]]$, $[v_{i_1},v_{i_2}]\in (\mathfrak{p})$.
Finally
$$[v_{i_1},\ldots,v_{i_m}]=-[v_{i_2},[[v_{i_3}\ldots,v_{i_m}],v_{i_1}]]
-[[v_{i_3},\ldots,v_{i_m}],[v_{i_1},v_{i_2}]]\in (\mathfrak{p}).$$

{\bf 2)} Let $i_1=i_2$. Applying the Jacobi identity to $[v_{i_1},v]$ we
get
$$\begin{array}{l}[v_{i_1},[v_{i_1},[v_{i_3},\overline{v}]]]=
-[v_{i_1},[v_{i_3},[\overline{v},v_{i_1}]]]-
[v_{i_1},[\overline{v},[v_{i_1},v_{i_3}]]]= \\
=[v_{i_3},[[\overline{v},v_{i_1}],v_{i_1}]]+
[[\overline{v},v_{i_1}],[v_{i_1},v_{i_3}]]+[\overline{v},[[v_{i_1},v_{i_3}],v_{i_1}]]
+[[\overline{v},v_{i_1}],[v_{i_1},v_{i_3}]].\end{array}$$
Note that, we have $\langle
e_{i_1},e_{i_1}+e_{\overline{v}}\rangle_q =2+\langle
e_{i_1},e_{\overline{v}}\rangle_q  \geq 2+(-1)= 1$, and therefore by
the induction hypothesis $[v_{i_1},[v_{i_1},\overline{v}]]\in (\mathfrak{p})$. Moreover,
$[v_{i_1},[v_{i_1},v_{i_3}]]\in(\mathfrak{r_1})\subseteq (\mathfrak{p})$. Finally
$[v_{i_1},v]=2[[v_{i_1},v_{i_3}][v_{i_1},\overline{v}]]$.
If $[v_{i_1},v_{i_3}]\in(\mathfrak{p})$, then
$[v_{i_1},v]\in (\mathfrak{p})$. Assume that $[v_{i_1},v_{i_3}]\not\in(\mathfrak{p})$.
In this case $\langle
e_{i_1},e_{i_3}\rangle_q  =-1$ and $$\langle e_{i_1},e_{\overline{v}}\rangle_q
= \langle e_{i_1},e_{i_2}+e_{i_3}+e_{\overline{v}}\rangle_q  -\langle e_{i_1},e_{i_2}\rangle_q-\langle
e_{i_1},e_{i_3}\rangle_q \geq 1-2-(-1)=0.$$  Then
$[v_{i_1},v]=2[[v_{i_1},v_{i_3}][v_{i_1},\overline{v}]]\in(\mathfrak{r}_0)\subseteq(\mathfrak{p})$
and we are done.

{\bf 3)} Let $\langle i_1, i_2\rangle_q=-1$.
By Lemma \ref{lem:pierwiastki} (c),
$\langle e_{i_1},e_{i_3}+\ldots+e_{i_m}\rangle\in\{-1,0,1\}$,
because $q$ is positive definite.
On the other hand
$$1\leq\langle e_{i_1},e_{i_2}+\ldots+e_{i_m}\rangle_q=
-1+\langle e_{i_1},e_{i_3}+\ldots+e_{i_m}\rangle_q\leq 0.$$
This contradiction shows that the case {\bf 3)} does not hold.

This finishes the proof. \epv

\begin{cor}\punkt If $q$ is a~positive definite quadratic form {\rm (\ref{eq:form})}, then
$$L(q,\mathfrak{r})\cong L(q,\mathfrak{p}). $$\epv
\end{cor}

\subsection{The second step of reduction}

Let $i_1,\ldots,i_m\in\{1,\ldots,n\}$. Following \cite{barotrivera}, we call the sequence
$(i_1,\ldots,i_m)$
a~\textbf{chordless cycle} of the form $q$
(\ref{eq:form}), if the following conditions are satisfied:

\begin{itemize}
 \item the elements $i_1,\ldots,i_m$ are pairwise different,
 \item $a_{ij}=\langle e_{i_j},e_{i_k}\rangle_q\neq 0$ if and only if $|k-j|=1$  mod $m$.
\end{itemize}

Chordless cycles are playing an~important role in \cite{barotrivera}, where
Lie algebras associated with positive definite quadratic forms are investigated.

A~chordless cycle $(i_1,\ldots,i_m)$ is called \textbf{positive}, if
$\langle e_{i_1},e_{i_m}\rangle_q=1$ and $\langle e_{i_j},e_{i_k}\rangle_q=-1$
for all $j,k$ such that $\{j,k\}\neq\{1,m\}$ and $|j-k|=1$ mod $m$.

\begin{rem}\punkt
 {\rm Note that if $(i_1,\ldots,i_m)$ is a~chordless cycle, then $(i_1,\ldots,i_m)$
is a~simple cycle in the bigraph $B(q)$. Moreover, if the
chordless cycle $(i_1,\ldots,i_m)$ is positive, then the cycle
$(i_1,\ldots,i_m)$ in $B(q)$ has exactly one broken edge
$\xymatrix{i_1\ar@{--}[r] & i_m}$. }
\end{rem}

Let $\mathfrak{r}_2\subseteq L(q)$ be the set consisting of all elements
$[v_{i_1},\ldots,v_{i_m}]$ such that $(i_1,\ldots,i_m)$ is a~positive chordless cycle.

\begin{lem}\punkt
  $\mathfrak{r}_2\subseteq \mathfrak{p}$.
\label{lem:chordless}\end{lem}

\textbf{Proof.} Let $v=[v_{i_1},\ldots,v_{i_m}]\in \mathfrak{r}_2$.
From the definition of a~positive chordless cycle, it follows easily
that $[v_{i_k},\ldots,v_{i_m}]$ is a~root for any $k>1$. Moreover
$\langle e_{i_1},e_{i_2}+\ldots+e_{i_m}\rangle_q=0$, and therefore $v\in\mathfrak{p}$ \epv

Set \begin{equation}
 \mathfrak{j}=\mathfrak{r}_1\cup \mathfrak{r}_2.
\end{equation}

For all elements $x,y\in L(q)$ we write $x\equiv y$ if $x-y\in
(\mathfrak{j})$. Obviously $\equiv$ is an equivalence relation.

Before we prove that the ideals $(\mathfrak{p})$ and $(\mathfrak{j})$ of $L(q)$
are equal, we need to prove two technical lemmata.

\begin{lem}\punkt
Let $q$ be a~positive definite quadratic form {\rm (\ref{eq:form})}. Let $m\geq 3$ be an~integer,
let $v=[v_{i_2},\ldots,v_{i_m}]$ be a~root and let
$(\mathfrak{j})_{m-1}=(\mathfrak{p})_{m-1}$. Let $i_1\in \{1,\ldots,n\}$
be such that $[v_{i_1},v_{i_2},\ldots,v_{i_m}]\in \mathfrak{p}$, $\langle e_{i_1},e_{i_2}\rangle_q=-1$ and
$\langle e_{i_1},e_{i_2}+\ldots+e_{i_m} \rangle_q= 0$. Then $[v_{i_1},v]\in(\mathfrak{j})$ or
there exists $\varepsilon\in\{0,1\}$ such that $[v_{i_1},v]\equiv (-1)^\varepsilon[v_{i_1},[a,x]]$, where

\ya $a=[v_{i_k},\ldots,v_{i_2}]$, $x=[v_{i_{k+1}},\ldots, v_{i_m}]$, for some $2\leq k<m$,

\yb $\langle e_{i_1},e_{i_j}\rangle_q=0$, for all $j=3,\ldots,k$, and

\yc $\langle e_{i_1},e_{i_{k+1}}\rangle_q=1$.
\label{lem:reduction-step1}\end{lem}

\textbf{Proof.} Let $m\geq 3$ and let $v=[v_{i_2},\ldots,v_{i_m}]$ be a~root. Let $i_1\in \{1,\ldots,n\}$
be such that $[v_{i_1},v_{i_2},\ldots,v_{i_m}]\in \mathfrak{p}$, $\langle e_{i_1},e_{i_2}\rangle_q=-1$ and
$\langle e_{i_1},e_{i_2}+\ldots+e_{i_m} \rangle_q= 0$. Note that $[v_{i_1},v]\equiv [v_{i_1},[\overline{a},x]]$,
where $\overline{a}=[v_{i_l},\ldots,v_{i_2}]$, $x=[v_{i_{l+1}},\ldots, v_{i_m}]$ and
$\langle e_{i_1},e_{i_j}\rangle_q=0$, for all $j=3,\ldots,l$ (i.e. the conditions
(a), (b) are satisfied). Indeed, it is enough to set
$l=2$, $\overline{a}=x_{i_2}$ and $x=[v_{i_3},\ldots, v_{i_m}]$.

Fix $\overline{a}$ and $\overline{x}$ such that $[v_{i_1},v]\equiv [v_{i_1},[\overline{a},\overline{x}]]$,
where $\overline{a}=[v_{i_l},\ldots,v_{i_2}]$, $\overline{x}=[v_{i_{l+1}},\ldots, v_{i_m}]$ and
$\langle e_{i_1},e_{i_j}\rangle_q=0$, for all $j=3,\ldots,l$. As we noted above there exists at least
one such a~presentation of  $[v_{i_1},v]$.
Consider the element $i_{l+1}$.
By Lemma \ref{lem:pierwiastki},  $\langle e_{i_1},e_{i_{l+1}}\rangle_q\in\{-1,0,1,2\}$.
\begin{itemize}
 \item If $\langle e_{i_1},e_{i_{l+1}}\rangle_q=1$, then we set $k=l$ and note that $[v_{i_1},v]$ has the
required form, i.e. $[v_{i_1},v]\equiv [v_{i_1},[a,x]]$, where $a=[v_{i_k},\ldots,v_{i_2}]$, $x=[v_{i_{k+1}},\ldots v_{i_m}]$,
 $\langle e_{i_1},e_{i_j}\rangle_q=0$, for all $j=3,\ldots,k$, and
$\langle e_{i_1},e_{i_{k+1}}\rangle_q=1$.

 \item If $\langle e_{i_1},e_{i_{l+1}}\rangle_q=-1$, then
$$\begin{array}{lcl}\langle e_{i_1},e_{i_{l+2}}+\ldots+e_{i_m}\rangle_q&=&\langle e_{i_1},e_v\rangle_q-
\langle e_{i_1},e_{\ov{a}}\rangle_q-\langle e_{i_1},e_{i_{l+1}}\rangle_q\\ &=&0-(-1)-(-1)=2.
\end{array}
$$
Therefore, by Lemma \ref{lem:pierwiastki}, $m=l+2$ and $i_1=i_{l+2}$. Note that

\begin{itemize}
 \item $[v_{i_1},[v_{i_1},v_{i_{l+1}}]]\in (\mathfrak{r}_1)\subseteq (\mathfrak{j})$,
 \item $\langle e_{i_1},e_{i_1}+e_{\overline{a}}\rangle_q=1$, then
 $[v_{i_1},[v_{i_1},\ov{a}]]\in
(\mathfrak{p})_{(m-1)}=(\mathfrak{j})_{(m-1)}\subseteq(\mathfrak{j})$,
 \item $\langle e_{i_1},e_{i_{l+1}}+e_{\overline{a}}\rangle_q=-2$, therefore
$[v_{i_{l+1}},\overline{a}]\in(\mathfrak{p})_{(m-1)}=(\mathfrak{j})_{(m-1)}\subseteq(\mathfrak{j})$.
\end{itemize}

Then we have
$$\begin{array}{lcl}[v_{i_1},v] &\equiv& [v_{i_1},[\ov{a},[v_{i_{l+1}},v_{i_1}]]]  \\
  &=& -[\ov{a},[[v_{i_{l+1}},v_{i_1}],v_{i_1}]]-[[v_{i_{l+1}},v_{i_1}],[v_{i_1},\ov{a}]] \\
  &\equiv & [[v_{i_1},\ov{a}],[v_{i_{l+1}},v_{i_1}]] \\
  &=& -[v_{i_{l+1}},[v_{i_1},[v_{i_1},\ov{a}]]]-[v_{i_1},[[v_{i_1},\ov{a}],v_{i_{l+1}}]] \\
  &\equiv & [v_{i_1},[v_{i_{l+1}},[v_{i_1},\ov{a}]]] \\
  &=& -[v_{i_1},[v_{i_1},[\ov{a},v_{i_{l+1}}]]]-[v_{i_1},[\ov{a},[v_{i_{l+1}},v_{i_1}]]] \\
  &\equiv & -[v_{i_1},[\ov{a},[v_{i_{l+1}},v_{i_1}]]]\\
  &\equiv&-[v_{i_1},v].
  \end{array}
$$

Therefore $2\cdot [v_{i_1},v]\in(\mathfrak{j})$ and $[v_{i_1},v]\in(\mathfrak{j})$.

\item If $\langle e_{i_1},e_{i_{l+1}}\rangle_q=2$, then $i_1=i_{l+1}$ and
$[v_{i_1},v]\equiv [v_{i_1},[\overline{a},[v_{i_1},y]]]$, where $y=[v_{i_{l+2}},\ldots,v_{i_m}]$.
By the bilinearity of $\langle -,-\rangle_q$ and assumptions, we have
$\langle e_{i_1},e_y\rangle_q=-1$.
Moreover
\begin{itemize}
 \item $\langle e_{i_1},e_{i_1}+e_y\rangle_q= 1$, then $[v_{i_1},[v_{i_1},y]]
 \in (\mathfrak{p})_{(m-1)}=(\mathfrak{j})_{(m-1)}\subseteq(\mathfrak{j})$,
 \item $\langle e_{i_1},e_{i_1}+e_{\overline{a}}\rangle_q=1$, then
 $[v_{i_1},[v_{i_1},\ov{a}]]\in
(\mathfrak{p})_{(m-1)}=(\mathfrak{j})_{(m-1)}\subseteq(\mathfrak{j})$,
 \item $\langle e_{i_1},e_y+e_{\overline{a}}\rangle_q=-2$, therefore
$[v_{i_{l+1}},\overline{a}]\in(\mathfrak{p})_{(m-1)}=(\mathfrak{j})_{(m-1)}\subseteq(\mathfrak{j})$.
\end{itemize}

Similarly as above
we can prove that $[v_{i_1},v]\in(\mathfrak{j})$.

\item Let $\langle e_{i_1},e_{i_{l+1}}\rangle_q=0$ and $y=[v_{i_{l+2}},\ldots,v_{i_m}]$. Consider
\footnotesize
$$ \begin{array}{lcl}[v_{i_1},[\overline{a},[v_{i_{l+1}},y]]]&=&-[v_{i_1},[v_{i_{l+1}},[y,\overline{a}]]]-
[v_{i_1},[y,[\overline{a},v_{i_{l+1}}]]] \\ &=&
[v_{i_{l+1}},[[y,\overline{a}],v_{i_1}]]+[[y,\overline{a}],[v_{i_1},v_{i_{l+1}}]]-
[v_{i_1},[[v_{i_{l+1}},\overline{a}],y]]. \end{array}
$$\normalsize

Note that $[v_{i_1},v_{i_{l+1}}]\in(\mathfrak{r}_1)\subseteq(\mathfrak{j})$,
$\langle e_{i_1},e_{\overline{a}}+e_y\rangle_q=0$ and therefore
$[[y,\overline{a}],v_{i_1}]\in(\mathfrak{p})_{m-1}=(\mathfrak{j})_{m-1}\subseteq (\mathfrak{j})$. Then
$$ [v_{i_1},v]\equiv [v_{i_1},[\overline{a},[v_{i_{l+1}},y]]]\equiv -[v_{i_1},[[v_{i_{l+1}},\overline{a}],y]]. $$ We may set $\overline{a}:=[v_{i_{l+1}},\overline{a}]$,
$\overline{x}:=y$ and continue this procedure inductively. 
\end{itemize}

Note that there exists $k$ such that $\langle e_{i_1},e_{i_{k+1}}\rangle_q\geq 1$, because $\langle e_{i_1},e_{i_2}+\ldots+e_{i_m}\rangle_q=0$
and $\langle e_{i_1},e_{i_2}\rangle_q=-1$.
%\end{itemize}
Therefore continuing this procedure inductively, we prove that $[v_{i_1},v]\in(\mathfrak{j})$ or
$[v_{i_1},v]\equiv (-1)^\varepsilon[v_{i_1},[a,x]]$, where $a=[v_{i_k},\ldots,v_{i_2}]$, $x=[v_{i_{k+1}},\ldots v_{i_m}]$,
 $\langle e_{i_1},e_{i_j}\rangle_q=0$, for all $j=3,\ldots,k$, and
$\langle e_{i_1},e_{i_{k+1}}\rangle_q=1$. \vsp
\epv

\begin{lem}\punkt  Assume that $q$ is a~positive definite quadratic form {\rm (\ref{eq:form})}. 
Let $m\geq 3$ be an~integer, $v=[v_{i_2},\ldots,v_{i_m}]$ be a~root and 
$(\mathfrak{j})_{m-1}=(\mathfrak{p})_{m-1}$. Let $i_1\in \{1,\ldots,n\}$
be such that $[v_{i_1},v_{i_2},\ldots,v_{i_m}]\in \mathfrak{p}$, $\langle e_{i_1},e_{i_2}\rangle_q=-1$ and
$\langle e_{i_1},e_{i_2}+\ldots+e_{i_m} \rangle_q= 0$. Moreover assume that
$[v_{i_1},v]\equiv (-1)^\varepsilon[v_{i_1},[a,x]]$ and the conditions {\rm (a)-(c)}
of Lemma {\rm \ref{lem:reduction-step1}} are satisfied. Then  $[v_{i_1},v]\in(\mathfrak{j})$ or
$[v_{i_1},v]\equiv [v_{i_1},[a,x]]\equiv [v_{i_1},[a,[b,y]]]$, where

\yi $a=[v_{i_k},\ldots,v_{i_2}]$, $b=[v_{i_s},v_{i_{s-1}}\ldots,v_{i_{k+1}}]$, $y=[v_{i_{s+1}},\ldots, v_{i_m}]$,

\yii $\langle e_{i_1},e_{i_j}\rangle_q=0$, for all $j=3,\ldots,k$ and $j=k+2,\ldots,s$,

\yiii $\langle e_{i_1},e_{i_2}\rangle_q=-1$, $\langle e_{i_1},e_{i_{k+1}}\rangle_q=1$,

%\yiv  for any $j\in{\rm supp}(b)$ we have $\langle e_j,e_a\rangle_q=0$, in particular 
\yiv $\langle e_b,e_a\rangle_q=0$,

\yv if $s<m$ then $\langle e_{i_{s+1}},e_a\rangle_q=-1$ and $\langle e_{i_{s+1}},e_b\rangle_q=-1$.

 \label{lem:reduction-step2}\end{lem}

\textbf{Proof.} Note that $[v_{i_1},v]\equiv [v_{i_1},[a,x]]\equiv [v_{i_1},[a,[b,y]]]$, where
$b=v_{i_{k+1}}$, $y=[v_{i_{k+2}},\ldots,v_{i_m}]$ and the conditions (i), (ii), (iii) are
satisfied, if we put $s=k+1$.
We may assume that the condition (iv) is also satisfied. Indeed, it is enough to show
that $\langle e_{i_{k+1}},e_a\rangle_q=0$.
\begin{itemize}
 \item If $\langle e_{i_{k+1}},e_a\rangle_q=-1$,
then $\langle e_{i_{k+1}},e_a+e_y\rangle_q=-2$. Therefore
$[a,y]\in(\mathfrak{r})_{m-1}=(\mathfrak{p})_{m-1}=(\mathfrak{j})_{m-1}\subseteq (\mathfrak{j})$. Then
$$\begin{array}{lcl} [v_{i_1},[a,[v_{i_{k+1}},y]]] &=& -[v_{i_1},[v_{i_{k+1}},[y,a]]]
- [v_{i_1},[y,[a,v_{i_{k+1}}]]]\\ & \equiv & [v_{i_1},[[a,v_{i_{k+1}}],y]] \\
  &=& -[[a,v_{i_{k+1}}],[y,v_{i_1}]]-[y,[v_{i_1},[a,v_{i_{k+1}}]]].
  \end{array}$$
Since $\langle e_{i_1},e_y\rangle_q=0=\langle
e_{i_1},e_{i_{k+1}}+e_a\rangle_q$, then $[y,v_{i_1}]\in
(\mathfrak{p})_{m-1}=(\mathfrak{j})_{m-1}$ and
$[v_{i_1},[a,v_{i_{k+1}}]]\in(\mathfrak{p})_{m-1}=(\mathfrak{j})_{m-1}$.Therefore
$[v_{i_1},[a,[v_{i_{k+1}},y]]]\in (\mathfrak{j})$.

\item  If
$\langle e_{i_{k+1}},e_a\rangle_q=2$, then $a=v_{i_{k+1}}$. It is a~contradiction, because
$$\langle e_{i_1},e_{i_{k+1}}\rangle_q=1\neq -1=\langle e_{i_1},e_a\rangle_q.$$

\item  If
$\langle e_{i_{k+1}},e_a\rangle_q=1$, then $\langle e_{i_{k+1}},e_a+e_y\rangle_q=0$.
Therefore $[v_{i_{k+1}},[y,a]]$, $[a,v_{i_{k+1}}]\in(\mathfrak{p})_{m-1}=(\mathfrak{j})_{m-1}$
and
$$[v_{i_1},[a,[v_{i_{k+1}},y]]] = -[v_{i_1},[v_{i_{k+1}},[y,a]]]
- [v_{i_1},[y,[a,v_{i_{k+1}}]]] \in(\mathfrak{j}).$$

\end{itemize}
Finally, we can assume that $\langle e_{i_{k+1}},e_a\rangle_q=\langle e_b,e_a\rangle_q=0$
and the condition (iv) is satisfied. Therefore $[v_{i_1},[a,[b,y]]]\equiv [v_{i_1},[b,[a,y]]]$,
because $\langle e_b,e_a\rangle_q=0$.\vsp

If $k+1=m$, then we are done. Assume that $k+1<m$ and consider the element $i_{k+2}$.

\begin{enumerate}
 \item Let $\langle e_{i_{k+2}},e_b\rangle_q=-1$.
\begin{enumerate}
 \item If $\langle e_{i_{k+2}},e_a\rangle_q=-1$,
then we put $s=k+1$ and we are done.
 \item If $\langle e_{i_{k+2}},e_a\rangle_q\geq 0$, then $[v_{i_{k+2}},a]\in(\mathfrak{j})_{m-1}$.

If $m=k+2$, then
$$ [v_{i_1},v]\equiv -[v_{i_1},[b,[v_{i_{k+2}},a]]]-[v_{i_1},[v_{i_{k+2}},[a,b]]]\in (\mathfrak{j})$$ and we are done.

Assume that $m>k+2$. We can assume that $\langle e_{i_{k+1}},e_a\rangle_q=0$.
Indeed, since $[v_{i_{k+2}},a]\in(\mathfrak{j})_{m-1}$, we have
$$\begin{array}{lcl}[v_{i_1},v]&\equiv& [v_{i_1},[b,[a,[v_{i_{k+2}},v_{i_{k+3}},\ldots,v_{i_m}]]]]\\ &\equiv& [v_{i_1},[b,[v_{i_{k+2}},[a,[v_{i_{k+3}},\ldots,v_{i_m}]]]]].
  \end{array}
$$
If $\langle e_{i_{k+2}},e_a\rangle_q\geq 1$, then
$\langle e_{i_{k+2}},e_a+e_{i_{k+3}}+\ldots+e_{i_m}\rangle_q\geq 1-1=0$. It follows that
$$[v_{i_{k+2}},[a,[v_{i_{k+3}},\ldots,v_{i_m}]]]\in (\mathfrak{p})_{m-1}=(\mathfrak{j})_{m-1}.$$
Therefore we can assume that $\langle e_{i_{k+2}},e_a\rangle_q=0$.
Moreover
$$\begin{array}{lcl} [v_{i_1},v]&\equiv& [v_{i_1},[b,[v_{i_{k+2}},[a,[v_{i_{k+3}},\ldots,v_{i_m}]]]]] \\
&=& -[v_{i_1},[v_{i_{k+1}},[[a,[v_{i_{k+3}},\ldots,v_{i_m}]],b]]] \\ &&-
[v_{i_1},[[a,[v_{i_{k+3}},\ldots,v_{i_m}]],[b,v_{i_{k+2}}]]]\\
&\equiv& -[v_{i_1},[[v_{i_{k+2}},b],[a,[v_{i_{k+3}},\ldots,v_{i_m}]]]] ,\end{array}
$$
because $\langle e_{i_{k+2}},e_b\rangle_q=-1$ and $\langle e_{i_{k+2}},e_a+e_b+e_{i_{k+3}}+\ldots +e_{i_m}\rangle_q=-2$.
We can assume that $\langle e_{i_1},e_{i_{k+2}}\rangle_q=0$. Indeed, if
$\langle e_{i_1},e_{i_{k+2}}\rangle_q=-1$, then
$$\begin{array}{lcl} 0=\langle e_{i_1},e_v\rangle_q&=&\langle e_{i_1},e_a+e_b+e_{i_{k+2}}\rangle_q+\langle e_{i_1},e_{i_{k+3}}+\ldots+e_{i_m}\rangle_q\\
&=&-1+\langle e_{i_1},e_{i_{k+3}}+\ldots+e_{i_m}\rangle_q.\end{array}$$
It follows that
$\langle e_{i_1},e_{i_{k+3}}+\ldots+v_{i_m}\rangle_q=1$, $\langle e_{i_1},e_a+e_{i_{k+3}}+\ldots+v_{i_m}\rangle_q=0$ and
$[v_{i_1},[a,[v_{i_{k+3}},\ldots,v_{i_m}]]]\in (\mathfrak{p})_{m-1}=(\mathfrak{j})_{m-1}$.
Therefore 
$$\begin{array}{lcl}[v_{i_1},v]&\equiv& -[v_{i_1},[[v_{i_{k+2}},b],[a,[v_{i_{k+3}},\ldots,v_{i_m}]]]]\\
&\equiv& -[[v_{i_{k+2}},b],[v_{i_1},[a,[v_{i_{k+3}},\ldots,v_{i_m}]]]]\in (\mathfrak{j}),
  \end{array}
$$ because $\langle e_{i_1},e_{i_{k+2}}+e_b\rangle_q=0$.
If $\langle e_{i_1},e_{i_{k+2}}\rangle_q\geq 1$, then
$$\begin{array}{lcl} 0=\langle e_{i_1},e_v\rangle_q&=&\langle e_{i_1},e_a+e_b+e_{i_{k+2}}\rangle_q+\langle e_{i_1},e_{i_{k+3}}+\ldots+e_{i_m}\rangle_q \\ &\geq& 1+\langle e_{i_1},e_{i_{k+3}}+\ldots+e_{i_m}\rangle_q.
  \end{array}
$$
It follows that $\langle e_{i_1},e_{i_{k+3}}+\ldots+v_{i_m}\rangle_q\leq -1$,
$\langle e_{i_1},e_a+e_{i_{k+3}}+\ldots+v_{i_m}\rangle_q\leq -2$ and $[a,[v_{i_{k+3}},\ldots,v_{i_m}]]\in (\mathfrak{p})_{m-1}=(\mathfrak{j})_{m-1}$.
Therefore, $[v_{i_1},v]\in (\mathfrak{j})$ and we are done. Finally,
$\langle e_{i_{k+2}},e_a\rangle_q=0$, $\langle e_{i_1},e_{i_{k+2}}\rangle_q=0$, $\langle e_{i_{k+2}}+e_b,e_a\rangle_q=0$ and
$$\begin{array}{lcl}[v_{i_1},v]&\equiv& -[v_{i_1},[[v_{i_{k+2}},b],[a,[v_{i_{k+3}},\ldots,v_{i_m}]]]] \\ &\equiv& -[v_{i_1},[a,[[v_{i_{k+2}},b],[v_{i_{k+3}},\ldots,v_{i_m}]]]].
  \end{array}
$$
We set $\ov{a}=a$, $\ov{b}=[v_{i_{k+2}},b]$, $\ov{y}=[v_{i_{k+3}},\ldots,v_{i_m}]$
and continue this procedure inductively using $[v_{i_1},[\overline{a},[\overline{b},\overline{y}]]]$ instead of $[v_{i_1},[a,[b,y]]]$.
\end{enumerate}

\item Let $\langle e_{i_{k+2}},e_b\rangle_q\neq -1$, then
$[v_{i_{k+2}},b]\in (\mathfrak{j})_{m-1}$. If $m=k+2$, then
$$[v_{i_1},v]\equiv [v_{i_1},[a,[b,v_{i_{k+2}}]]]\in(\mathfrak{j}) $$ and we are done.
Assume that $m>k+2$. Since $[v_{i_{k+2}},b]\in (\mathfrak{j})_{m-1}$, we have
$$ [v_{i_1},[a,[b,[v_{i_{k+2}},[v_{i_{k+3}},\ldots,v_{i_m}]]]]] \equiv
[v_{i_1},[a,[v_{i_{k+2}},[b,[v_{i_{k+3}},\ldots,v_{i_m}]]]]].$$
For the sake of simplicity we present partial results in tables.
In the first column of the~following table we consider all possible values of
$\langle e_{i_{k+2}},e_b\rangle_q$. In the second column we give the corresponding value
of $\langle e_{i_{k+2}},e_b+e_{i_{k+3}}+\ldots +e_{i_m}\rangle_q$. The third column contains
the sign "+", if we can deduce that
$X=[v_{i_{k+2}},[b,[v_{i_{k+3}},\ldots,v_{i_m}]]]\in(\mathfrak{j})_{m-1}$, and the sign "$-$", otherwise.\vsp

\noindent \begin{tabular}{|ccccccc|}
\hline
& $\langle e_{i_{k+2}},e_b\rangle_q$ &\vline & $\langle e_{i_{k+2}},e_b+e_{i_{k+3}}+\ldots +e_{i_m}\rangle_q$ &\vline&   $X\in (\mathfrak{j})_{m-1}$ &\\
\hline & 0 &\vline& -1 &\vline& $-$ & \\
\hline & 1 &\vline&  0 &\vline& $+$ & \\
\hline & 2 &\vline&  1 &\vline& $+$ &  \\
\hline
\end{tabular} \vsp

\noindent Therefore we may assume that $\langle e_{i_{k+2}},e_b\rangle_q=0$, because
otherwise $$[v_{i_1},v]\equiv [v_{i_1},[a,[v_{i_{k+2}},[b,[v_{i_{k+3}},\ldots,v_{i_m}]]]]]\in(\mathfrak{j}).$$

\begin{enumerate}
 \item Assume that $\langle e_{i_{k+2}},e_a\rangle_q\neq -1$. Then
$[v_{i_{k+2}},a]\in (\mathfrak{j})_{m-1}$ and
$$\begin{array}{lcl} [v_{i_1},v] &\equiv&
[v_{i_1},[a,[v_{i_{k+2}},[b,[v_{i_{k+3}},\ldots,v_{i_m}]]]]] \\ &\equiv& [v_{i_1},[v_{i_{k+2}},[a,[b,[v_{i_{k+3}},\ldots,v_{i_m}]]]]].
\end{array}
$$
In the first column of the~following table we consider all possible values of
$\langle e_{i_{k+2}},e_a\rangle_q$. In the second column we give the corresponding value
of $x=\langle e_{i_{k+2}},e_a+e_b+e_{i_{k+3}}+\ldots +e_{i_m}\rangle_q$. The third column contains
the sign "+", if we can deduce that
$$X=[v_{i_{k+2}},[a,[b,[v_{i_{k+3}},\ldots,v_{i_m}]]]]\in(\mathfrak{j})_{m-1},$$ and the sign 
"$-$", otherwise. \vsp

\noindent \begin{tabular}{|ccccccc|}
\hline
& $\langle e_{i_{k+2}},e_a\rangle_q$ &\vline & $x$ &\vline&  $X\in (\mathfrak{j})_{m-1}$ &\\
\hline & 0 &\vline& -1 &\vline& $-$ & \\
\hline & 1 &\vline&  0 &\vline& $+$ & \\
\hline & 2 &\vline&  1 &\vline& $+$ &  \\
\hline
\end{tabular} \vsp

\noindent Therefore we may assume that $\langle e_{i_{k+2}},e_a\rangle_q=0$, because
otherwise $[v_{i_1},v]\in(\mathfrak{j})$.
Moreover,\small
$$[v_{i_1},[v_{i_{k+2}},[a,[b,[v_{i_{k+3}},\ldots,v_{i_m}]]]]]\equiv
[v_{i_1},[v_{i_{k+2}},[b,[a,[v_{i_{k+3}},\ldots,v_{i_m}]]]]], $$ \normalsize
because $\langle e_a,e_b\rangle_q=0$.

In the first column of the~following table we consider all
possible values of $\langle e_{i_1},e_{i_{k+2}}\rangle_q$. In the
second column we present a~consequences of the information
contained in the first column. Finally, in the second table we
present conclusions of the results presented in the first table.\vsp

\noindent \begin{tabular}{|cccc|}
\hline
& $\langle e_{i_1},e_{i_{k+2}}\rangle_q$ &\vline & consequences \\
\hline & -1 &\vline& $\langle e_{i_1},e_b+e_{i_{k+3}}+\ldots +e_{i_m}\rangle_q=2$ \\
\hline & 0 &\vline& $\langle e_{i_1},e_a+e_b+e_{i_{k+3}}+\ldots +e_{i_m}\rangle_q=0$  \\
\hline & 1 &\vline& $\langle e_{i_1},e_a+e_{i_{k+3}}+\ldots +e_{i_m}\rangle_q=-2$   \\
\hline & 2 &\vline& $\langle e_{i_1},e_a+e_b+e_{i_{k+3}}+\ldots +e_{i_m}\rangle_q=-2$    \\
\hline
%\end{tabular}
%
%\noindent \begin{tabular}{|cccc|}
\hline
& $\langle e_{i_1},e_{i_{k+2}}\rangle_q$ &\vline & conclusions  \\
\hline & -1 &\vline&
$[b,[v_{i_{k+3}},\ldots,v_{i_m}]]\in(\mathfrak{j})_{m-1}$  \\
\hline & 0 &\vline&
$[v_{i_1},[a,[b,[v_{i_{k+3}},\ldots,v_{i_m}]]]]\in(\mathfrak{j})_{m-1}$  \\
\hline & 1 &\vline&  $[a,[v_{i_{k+3}},\ldots,v_{i_m}]]\in(\mathfrak{j})_{m-1}$  \\
\hline & 2 &\vline&  $[a,[b,[v_{i_{k+3}},\ldots,v_{i_m}]]]\in(\mathfrak{j})_{m-1}$   \\
\hline
\end{tabular} \vsp

All these cases imply that $[v_{i_1},v]\in(\mathfrak{j})$.

\item Assume that $\langle e_{i_{k+2}},e_a\rangle_q=-1$.
In this case
$$\langle e_{i_{k+2}},e_a+e_b+e_{i_{k+3}}+\ldots+e_{i_m}\rangle_q=-2.$$ It follows that
$$[a,[b,[v_{i_{k+3}},\ldots,v_{i_m}]]]\in(\mathfrak{j})_{m-1}$$ and

$$\begin{array}{lcl}[v_{i_1},v] &\equiv &
[v_{i_1},[a,[v_{i_{k+2}},[b,[v_{i_{k+3}},\ldots,v_{i_m}]]]]] \\
& = & [v_{i_1},[v_{i_{k+2}},[a,[b,[v_{i_{k+3}},\ldots,v_{i_m}]]]]] \\
&&- [v_{i_1},[[b,[v_{i_{k+3}},\ldots,v_{i_m}]],[a,v_{i_{k+2}}]]] \\
& \equiv & -[v_{i_1},[[v_{i_{k+2}},a],[b,[v_{i_{k+3}},\ldots,v_{i_m}]]]].
\end{array}
$$
Consider $\langle e_{i_1}, e_{i_{k+2}}\rangle_q$ and \begin{footnotesize}
$$\begin{array}{l}
 -[v_{i_1},[[v_{i_{k+2}},a],[b,[v_{i_{k+3}},\ldots,v_{i_m}]]]]= \\
= [[v_{i_{k+2}},a],[[b,[v_{i_{k+3}},\ldots,v_{i_m}]],v_{i_1}]]
+[[b,[v_{i_{k+3}},\ldots,v_{i_m}]],[v_{i_1},[v_{i_{k+2}},a]]].\end{array}$$
\end{footnotesize}
We present again partial results in tables.
\vsp

\noindent \begin{tabular}{|cccc|}
\hline
& $\langle e_{i_1},e_{i_{k+2}}\rangle_q$ &\vline & consequences \\
\hline & -1 &\vline& $\langle e_{i_1},e_b+e_{i_{k+3}}+\ldots +e_{i_m}\rangle_q=2$  \\
\hline & 1 &\vline&
\begin{tabular}{c} $\langle e_{i_1},e_b+e_{i_{k+3}}+\ldots +e_{i_m}\rangle_q=0$ \\
and $\langle e_{i_1},e_{i_{k+2}}+e_a\rangle_q=0$ \end{tabular}
   \\
\hline & 2 &\vline& $i_1=i_{k+2}$   \\
\hline
%\end{tabular}
%
%\noindent \begin{tabular}{|cccc|}
\hline
& $\langle e_{i_1},e_{i_{k+2}}\rangle_q$ &\vline& conclusions  \\
\hline & -1 &\vline&
$[b,[v_{i_{k+3}},\ldots,v_{i_m}]]\in(\mathfrak{j})_{m-1}$   \\
\hline & 1 &\vline&
 \begin{tabular}{c}  $[[b,[v_{i_{k+3}},\ldots,v_{i_m}]],v_{i_1}]\in(\mathfrak{j})_{m-1}$\\
and $[v_{i_1},[v_{i_{k+2}},a]]\in(\mathfrak{j})_{m-1}$
\end{tabular}  \\
\hline & 2  &\vline& \begin{tabular}{c} contradiction, because \\
$\langle e_{i_1},e_b\rangle_q=1\neq 0=\langle e_{i_{k+2}},e_b\rangle_q$ \end{tabular}
  \\
\hline
\end{tabular} \vsp

\noindent Therefore, we can assume that $\langle e_{i_1},e_{i_{k+2}}\rangle_q=0$, because
otherwise $[v_{i_1},v]\in(\mathfrak{j})$. Moreover
$$ \begin{array}{lcl}
[v_{i_1},v]&\equiv& -[v_{i_1},[[v_{i_{k+2}},a],[b,[v_{i_{k+3}},\ldots,v_{i_m}]]]]\\
&\equiv& -[v_{i_1},[b,[[v_{i_{k+2}},a],[v_{i_{k+3}},\ldots,v_{i_m}]]]],
\end{array}$$
because $\langle e_{i_{k+2}}+e_a,e_b\rangle_q=0$. Therefore
$$[v_{i_1},v]\equiv -[v_{i_1},[\ov{a},[\ov{b},\ov{y}]]], $$
where $\ov{a}=[v_{i_{k+2}},a]$, $\ov{b}=b$, $\ov{y}=[v_{i_{k+3}},\ldots,v_{i_m}]$
and the conditions (i)-(iv) are satisfied.
\end{enumerate}
\end{enumerate}

Continuing this procedure inductively we show that $[v_{i_1},v]\in(\mathfrak{j})$ or
$[v_{i_1},v]\equiv [v_{i_1},[a,x]]\equiv [v_{i_1},[a,[b,y]]]$ and the conditions (i)-(v)
are satisfied.\epv

\begin{prop}\punkt
Let $q$ be a~positive definite quadratic form. The ideals $(\mathfrak{j})$ and $(\mathfrak{p})$  of the Lie algebra $L(q)$
are equal.
\label{prop:second-reduction}\end{prop}

\textbf{Proof.} The inclusion $(\mathfrak{j})\subseteq
(\mathfrak{p})$ is obvious. It is enough to prove that
$(\mathfrak{p})\subseteq (\mathfrak{j})$. Let
$v=[v_{i_2},\ldots,v_{i_m}]$ be a~root in $L(q)$ and let $i_1\in
\{1,\ldots,n\}$ be such that $[v_{i_1},v_{i_2},\ldots,v_{i_m}]\in
\mathfrak{p}$ and $\langle e_{i_1},e_{i_2}+\ldots+e_{i_m}
\rangle_q= 0$. We claim that $[v_{i_1},v]\in \mathfrak{j}$.

We prove our claim by induction on $\ell(v)=m-1$. For $\ell(v)< 3$ our
statement easily follows by a~case by case inspection on all
possible cases.

Let $\ell(v)\geq 3$, then $m\geq 4$, $v=[v_{i_2},[v_{i_3},
\overline{v}]]$, $\langle e_{i_2},e_{i_3}+e_{\overline{v}}\rangle_q  =-1$,
$\langle e_{i_3},e_{\overline{v}}\rangle_q  =-1$ and
$\ell(\overline{v})\geq 1$. By Lemma \ref{lem:pierwiastki}, it is enough to
consider the following three cases. \vsp

\textbf{1)} If $\langle e_{i_1},e_{i_2}\rangle_q=0$, then by the bilinearity of
$\langle -,-\rangle_q$, we have $\langle e_{i_1},e_{i_3}+\ldots+e_{i_m}\rangle_q=0$.
Moreover $$[v_{i_1},v]=-[v_{i_2},[[v_{i_3},\ldots,v_{i_m}],v_{i_1}]]-[[v_{i_3},
\ldots,v_{i_m}],[v_{i_1},v_{i_2}]]. $$
By definitions, $[v_{i_1},v_{i_2}]\in \mathfrak{j}$ and
$[v_{i_1},[v_{i_3},\ldots,v_{i_m}]]\in \mathfrak{r}$. Then, by Proposition \ref{prop:first-red},
$[v_{i_1},[v_{i_3},\ldots,v_{i_m}]]\in (\mathfrak{p})$ and by the induction hypothesis
$[v_{i_1},[v_{i_3},\ldots,v_{i_m}]]\in (\mathfrak{j})$. Finally, $[v_{i_1},v]\in (\mathfrak{j})$.\vsp

\textbf{2)} If $i_1=i_2$, then $\langle e_{i_1},e_{i_2}\rangle_q=2$ and
$\langle e_{i_1},e_{i_3}+\ldots+e_{i_m}\rangle_q=-2$. This is a~contradiction with Lemma
\ref{lem:pierwiastki} and therefore the case 2) does not hold.\vsp

\textbf{3)} Let $\langle e_{i_1},e_{i_2}\rangle_q\in\{ 1,-1\}$. In this case
we apply the Jacobi identity and develop Lemmata \ref{lem:reduction-step1}, \ref{lem:reduction-step2}
to find an element $w\in L(q)$ such that
$[v_{i_1},v]-w\in(\mathfrak{j})$ 
(i.e. $[v_{i_1},v]\equiv w$). Finally we show that $w\in (\mathfrak{j})$, which implies
that $[v_{i_1},v]\in(\mathfrak{j})$.
 \vsp

\textbf{3.1)} Let $\langle e_{i_1},e_{i_2}\rangle_q=1$. We reduce this case to the case
3.2) presented below. Since
$\langle e_{i_1},e_v\rangle_q=0$ and $\langle e_{i_1},e_{i_2}\rangle_q=1$, then there exists
$k\in\{3,\ldots,m\}$ such that $\langle
e_{i_1},e_{i_k}\rangle_q=-1$. Choose $k$ minimal with this
property. We may assume that $k<m$, because $[v_{i_{m-1}},v_{i_m}]=-[v_{i_m},v_{i_{m-1}}]$
and we can work with $-v$ instead of $v$.
Note that for all
$s=3,\ldots,k-1$, we have $\langle e_{i_1},e_{i_s}\rangle_q=0$.
Indeed, if there exists $s=3,\ldots,k-1$ such that $\langle
e_{i_1},e_{i_s}\rangle_q\neq 0$, then by the choice of $k$, we
have $\langle e_{i_1},e_{i_s}\rangle_q\geq 1$. Then $\langle
e_{i_1},e_{i_k}+\ldots+e_{i_m}\rangle_q\leq \langle
e_{i_1},e_v\rangle_q-\langle e_{i_1},e_{i_2}\rangle_q-\langle
e_{i_1},e_{i_s}\rangle_q=-2$ and we get a~contradiction, because
$[v_{i_k},\ldots,v_{i_m}]$ is a~root.

Now applying the Jacobi identity we get
$$[v_{i_2},\ldots,v_{i_{k-1}},[v_{i_k},y]]=[v_{i_2},\ldots,v_{i_k},[v_{i_{k-1}},y]]
+[v_{i_2},\ldots,v_{i_{k-2}},[[v_{i_{k-1}},v_{i_k}],y]]. $$
By Lemma \ref{lem:root} (a), $[v_{i_k},[v_{i_{k-1}},y]]$ or $[v_{i_{k-1}},v_{i_k}]$ is not a~root,
then $[v_{i_{k-1}},v_{i_k}]\in(\mathfrak{p})$ or $[v_{i_k},[v_{i_{k-1}},y]]\in(\mathfrak{p})$.
By the induction hypothesis $[v_{i_k},[v_{i_{k-1}},y]]\in(\mathfrak{j})$ or $[v_{i_{k-1}},v_{i_k}]\in(\mathfrak{j})$.
Therefore $v\equiv [v_{i_2},\ldots,v_{i_{k-2}},[x,z]]$, where $x=v_{i_k}$ and $z=[v_{i_{k-1}},y]$
or $x=[v_{i_{k-1}},v_{i_k}]$, $z=y$.
In both cases $\langle e_{i_1},e_x\rangle_q=-1$.
Continuing this procedure (i.e. $x$ plays a~role of $v_{i_k}$ and $z$ plays a~role of $y$),
we get
$$[v_{i_1},v]\equiv [v_{i_1},[v_{i_2},[x,z]]], $$ where $\langle e_{i_1},e_x\rangle_q=-1$.
Applying the Jacobi identity, we get
$$[v_{i_2},[x,z]]=[x,[v_{i_2},z]]+[[v_{i_2},x],z]. $$ By Lemma
\ref{lem:root} (b), $[v_{i_2},z]$ or $[v_{i_2},x]$ is not a~root,
then $[x,[v_{i_2},z]]\in(\mathfrak{p})$ or
$[v_{i_2},x]\in(\mathfrak{p})$. By the induction hypothesis
$[x,[v_{i_2},z]]\in(\mathfrak{j})$ or $[v_{i_2},x]\in(\mathfrak{j})$.
Therefore $[v_{i_1},v]\equiv[v_{i_1},[x,[v_{i_2},z]]]$ or
$[v_{i_1},v]\equiv [v_{i_1},[[v_{i_2},x],z]]$. If
$[v_{i_1},v]\equiv[v_{i_1},[x,[v_{i_2},z]]]$, then applying Lemma
\ref{lem:root1} we get a~reduction to the case 1) or to the case
3.2) below. If $[v_{i_1},v]\equiv [v_{i_1},[[v_{i_2},x],z]]$, then
$$[v_{i_1},v]\equiv [v_{i_1},[[v_{i_2},x],z]]=-[[v_{i_2},x],[z,v_{i_1}]]-[z,[v_{i_1},[v_{i_2},x]]]. $$
Note that $\ell(z)\geq 1$, because we choose $k$ with the property $k<m$. Then
$\langle e_{i_1},e_z\rangle_q=\langle e_{i_1},e_v\rangle_q-\langle e_{i_1},e_x\rangle_q-\langle e_{i_1},e_{i_2}\rangle_q=0-(-1)-1=0$ and $\langle e_{i_1},e_{i_2}+e_x\rangle_q=1-1=0$, and therefore
by the induction hypothesis 
$$[v_{i_1},v]=-[[v_{i_2},x],[z,v_{i_1}]]-[z,[v_{i_1},[v_{i_2},x]]]\in(\mathfrak{j}).$$

\textbf{3.2)} Let $\langle e_{i_1},e_{i_2}\rangle_q=-1$. Applying Lemmata \ref{lem:reduction-step1}
and \ref{lem:reduction-step2} we get  $[v_{i_1},v]\in(\mathfrak{j})$ or
$[v_{i_1},v]\equiv [v_{i_1},[a,x]]\equiv [v_{i_1},[a,[b,y]]]$, where

(i) $a=[v_{i_k},\ldots,v_{i_2}]$, $b=[v_{i_s},v_{i_{s-1}}\ldots,v_{i_{k+1}}]$, $y=[v_{i_{s+1}},\ldots v_{i_m}]$,

(ii) $\langle e_{i_1},e_{i_j}\rangle_q=0$, for all $j=3,\ldots,k$ and $j=k+2,\ldots,s$,

(iii) $\langle e_{i_1},e_{i_2}\rangle_q=-1$, $\langle e_{i_1},e_{i_{k+1}}\rangle_q=1$,

(iv) $\langle e_b,e_a\rangle_q=0$,

(v) if $s<m$ then $\langle e_{i_{s+1}},e_a\rangle_q=-1$ and $\langle e_{i_{s+1}},e_b\rangle_q=-1$. \vsp

Consider the following cases.

\textbf{(a)} If $s=m$, then $[v_{i_1},v]\equiv [v_{i_1},[a,b]]$ and
$$q(e_a+e_b)=q(e_a)+q(e_b)+\langle e_a,e_b\rangle_q=1+1+0=2.$$
Therefore $e_a+e_b$ is not a~root of $q$, by the induction hypothesis $[a,b]\in(\mathfrak{j})$
and $[v_{i_1},v]\equiv[v_{i_1},[a,b]]\in(\mathfrak{j})$.\vsp

We may assume that $s<m$ and consider $\langle e_{i_1},e_{i_{s+1}}\rangle_q$.
Partial results are presented in the following tables.\vsp

\noindent \begin{tabular}{|cccc|}
\hline
& $\langle e_{i_1},e_{i_{s+1}}\rangle_q$ &\vline & consequence \\
\hline & -1 &\vline& $\begin{array}{l} q(e_{i_1}+(e_{i_{s+1}}+e_a)) =\\
q(e_{i_1})+q(e_{i_{s+1}}+e_a)+\langle e_{i_1},e_{i_{s+1}}+e_a\rangle_q =\\
1+1-1-1=0\end{array}$   \\
\hline & 1 &\vline&
$\begin{array}{l} q(-e_{i_1}+(e_{i_{s+1}}+e_b)) =\\
q(e_{i_1})+q(e_{i_{s+1}}+e_b)-\langle e_{i_1},e_{i_{s+1}}+e_b\rangle_q =\\
1+1-1-1=0\end{array}$
  \\
\hline & 2 &\vline& $i_1=i_{s+1}$   \\
\hline
%\end{tabular}
%
%\noindent \begin{tabular}{|cccc|}
\hline
& $\langle e_{i_1},e_{i_{s+1}}\rangle_q$ &\vline & conclusions  \\
\hline & -1 &\vline&  \begin{tabular}{c}
contradiction, because \\ $q$ is positive definite \end{tabular}   \\
\hline & 1 &\vline& \begin{tabular}{c}
contradiction, because \\ $q$ is positive definite \end{tabular}  \\
\hline & 2  &\vline& \begin{tabular}{c} contradiction, because \\
$\langle e_{i_{s+1}},e_b\rangle_q=-1\neq 1=\langle e_{i_1},e_b\rangle_q$ \end{tabular}  \\
\hline
\end{tabular} \vsp

Therefore, we may assume that $\langle e_{i_1},e_{i_{s+1}}\rangle_q=0$.

\textbf{(b)} Let $m\geq s+2$. Then

$$[v_{i_1},v]\equiv [v_{i_1},[a,[b,[v_{i_{s+1}},z]]]], $$
where $\ell(z)\geq 1$. Moreover $\langle e_{i_{s+1}},e_b+e_z\rangle_q=-1-1=-2$, then
$[z,b]$ is not a~root and by the induction hypothesis $[z,b]\in\mathfrak({j})$. Therefore
$$ \begin{array}{lcl}[v_{i_1},[a,[b,[v_{i_{s+1}},z]]]]
&=& -[v_{i_1},[a,[v_{i_{s+1}},[z,b]]]]-[v_{i_1},[a,[z,[b,v_{i_{s+1}}]]]]\\ 
&\equiv& [v_{i_1},[a,[[b,v_{i_{s+1}}],z]]].
   \end{array}
$$
Applying the Jacobi identity, we get
$$[v_{i_1},[a,[[b,v_{i_{s+1}}],z]]]=-[v_{i_1},[[b,v_{i_{s+1}}],[z,a]]]-
[v_{i_1},[z,[a,[b,v_{i_{s+1}}]]]]. $$
Note that $\langle e_{i_{s+1}},e_a+e_z\rangle_q=-1-1=-2$, then $[z,a]$
is not a~root and by the induction hypothesis $[z,a]\in(\mathfrak{j})$. Therefore
$$\begin{array}{lcl}
[v_{i_1},[a,[[b,v_{i_{s+1}}],z]]]&\equiv& [v_{i_1},[[a,[b,v_{i_{s+1}}]],z]] \\ &=&
-[[a,[b,v_{i_{s+1}}]],[z,v_{i_1}]]-[z,[v_{i_1},[a,[b,v_{i_{s+1}}]]]].
\end{array}
 $$
Note that $\langle e_{i_1},e_a+e_b+e_{i_{s+1}}\rangle_q=-1+1+0=0$ and
$\langle e_{i_1},e_z\rangle_q=\langle e_{i_1},e_v\rangle_q-\langle e_{i_1},e_a+e_b+e_{i_{s+1}}\rangle_q=0$.
By the induction hypothesis $[v_{i_1},[a,[b,v_{i_{s+1}}]]]\in(\mathfrak{j})$ and
$[z,v_{i_1}]\in(\mathfrak{j})$. Therefore $[v_{i_1},v]\in(\mathfrak{j})$.\vsp

\textbf{(c)} Let $m=s+1$. We recall that $\langle e_{i_1},e_{i_2}\rangle_q=-1$,
$\langle e_{i_1},e_{i_{s+1}}\rangle_q=1$ and $\langle e_{i_1},e_{i_j}\rangle_q=0$
for all $j=3,\ldots,s$. Applying 1) and
the Jacobi identity, it is straightforward to
prove the following conditions:

(i) $ [v_{i_1},v]\equiv(-1)^\varepsilon 
[v_{i_1},v_{i_2},\ldots,v_{i_k},v_{i_{s+1}},v_{i_s},\ldots,v_{i_{k+1}}]$,

(ii) $[v_{i_1},v]\equiv(-1)^\varepsilon
[v_{i_1},[[v_{i_j},\ldots,v_{i_2}],v_{i_{j+1}},\ldots,v_{i_{k+1}}]]$, for
 all
 $j=2,\ldots,k,s+1$, $s\ldots,k+2$.

Without loss of generality, we can assume that in both cases $\varepsilon=0$.
It follows from (i) that there exists a~numbering of elements $\{i_2,\ldots,i_{s+1}\}$,
such that $$[v_{i_1},v]\equiv [v_{i_1},v_{i_2},\ldots,v_{i_{s+1}}], $$
where, $[v_{i_2},\ldots,v_{i_{s+1}}]$ is a~root, $\langle e_{i_1},e_{i_2}\rangle_q=-1$,
$\langle e_{i_1},e_{i_{s+1}}\rangle_q=1$ and $\langle e_{i_1},e_{i_j}\rangle_q=0$
for $j=3,\ldots,s$. Moreover, it follows from (ii) that the elements
$[v_{i_j},\ldots,v_{i_2}]$ for all  $j=2,\ldots,k,s+1$, $s\ldots,k+2$, are roots,
because otherwise $[v_{i_j},\ldots,v_{i_2}]\in\mathfrak{j}$ and $[v_{i_1},v]\in\mathfrak{j}$.

We claim that $\langle e_{i_j},e_{i_{j+1}}\rangle_q=-1$ for all $j=2,\ldots,s$.
Assume, for the contrary,
that there exists $j=2,\ldots,s$, such that $\langle e_{i_j},e_{i_{j+1}}\rangle\neq -1$.
If $j=s$, then $[v_{i_2},\ldots,v_{i_{s+1}}]$ is not a~root. If $j=2$, then
$[v_{i_2},v_{i_3}]\in\mathfrak{j}$. Applying the Jacobi identity we get
$$[v_{i_1},[v_{i_2},\ldots,v_{i_{s+1}}]]\equiv
[v_{i_1},[v_{i_3},[v_{i_2},[v_{i_4},\ldots,v_{i_{s+1}}]]]] $$
and, by the case 1), $[v_{i_1},v]\in\mathfrak{j}$. Therefore we can assume that $j=3,\ldots,s-1$.
Since $[v_{i_2},\ldots,v_{i_{s+1}}]$ is a~root,
$\langle e_{i_j},e_{i_{j+1}}+\ldots +e_{i_{s+1}}\rangle_q=-1$ and
$\langle e_{i_j},e_{i_{j+2}}+\ldots +e_{i_{s+1}}\rangle_q\geq-1$. By the bilinearity of
$\langle -,-\rangle_q$, Lemma \ref{lem:pierwiastki}(c) and our assumptions, we have $\langle e_{i_j},e_{i_{j+2}}+\ldots +e_{i_{s+1}}\rangle=-1$
and $\langle e_{i_j},e_{i_{j+1}}\rangle=0$. By assumptions and (ii), the elements
$[v_{i_j},v_{i_{j-1}},\ldots,v_{i_2}]$ and $[v_{i_{j+1}},v_{i_{j+2}},\ldots,v_{i_{s+1}}]$ are roots.
Moreover $\langle e_{i_{j+1}},e_{i_2}+\ldots +e_{i_{j-1}}\rangle_q=\langle e_{i_{j+1}},e_{i_2}+\ldots +e_{i_{j-1}}+e_{i_j}\rangle_q=-1$.
Set $y=[v_{i_{j-1}},\ldots,v_{i_2}]$ and $x=[v_{i_{j+2}},\ldots,v_{i_{s+1}}]$, then
$$
\begin{array}{lcl}
 1=q(v)&=&q(e_y+e_{i_j}+e_{i_{j+1}}+e_x) \\
  &=& q(e_y+e_{i_j})+q(e_{i_{j+1}}+e_x)+\langle e_y+e_{i_j},e_{i_{j+1}}+e_x\rangle_q \\
  &=& 2+\langle e_y,e_{i_{j+1}}\rangle_q+\langle e_y,e_x\rangle_q+
         \langle e_{i_j},e_{i_{j+1}}\rangle_q+\langle e_{i_j},e_x\rangle_q \\
  &=& 2+ (-1)+ \langle e_y,e_x\rangle_q +0+(-1) \\
  &=& \langle e_a,e_x\rangle_q.
\end{array}
$$
It follows
$$
\begin{array}{lcl}
  q(-e_x+(e_y+e_{i_1}))&=&q(e_x)+q(e_y+e_{i_1})-\langle e_x,e_y+e_{i_1}\rangle_q\\
 &=& 1+1-\langle e_x,e_y\rangle_q -\langle e_x,e_{i_1}\rangle_q\\
&=&2-1-1=0,
\end{array}
$$
because $\langle e_x,e_{i_1}\rangle_q=\langle e_{i_{s+1}},e_{i_1}\rangle_q=1$.
This is a~contradiction, because $q$ is positive definite.
Finally, we proved that $\langle e_{i_j},e_{i_{j+1}}\rangle_q=-1$ for all $j=2,\ldots,s$.

If $\langle e_{i_j},e_{i_l}\rangle_q\leq 0$ for all $2\leq j<l \leq s+1$, then
$(i_1,i_2,\ldots,i_{s+1})$ is a~positive chordless cycle and therefore $[v_{i_1},v_{i_2},\ldots,v_{i_{s+1}}]\in\mathfrak{j}$. Indeed, if
$(i_1,i_2,\ldots,i_{s+1})$ is not a~positive chordless cycle, then there
exists $2\leq j<l\leq s+1$ such that $l\neq j,j+1$ and $\langle e_{i_j},e_{i_l}\rangle_q=-1$. 
Therefore $q(e_{i_2}+\ldots+e_{i_{s+1}})\leq s-(s-1)+\langle e_{i_j},e_{i_l}\rangle_q=0$ and 
$q$ is not positive definite.

Assume that $\langle e_{i_j},e_{i_l}\rangle_q> 0$ for some $2\leq j<l \leq s+1$.
%Note that there exists $2\leq j<l \leq s+1$ such that $\langle e_{i_j},e_{i_l}\rangle_q=1$.
%Indeed, if $\langle e_{i_j},e_{i_l}\rangle_q\neq 1$ for all $2\leq j<l \leq s+1$, then
%there exists $2\leq j<l \leq s+1$ such that $\langle e_{i_j},e_{i_l}\rangle_q=2$. 
Choose $j,l$ such that $2\leq j<l-1 \leq s+1$ and $l-j$ is minimal with the property
$\langle e_{i_j},e_{i_l}\rangle_q\neq 0$. 

If $\langle e_{i_j},e_{i_l}\rangle_q=-1$, then $q(e_{i_j}+\ldots+e_{i_l})=0$.
If $\langle e_{i_j},e_{i_l}\rangle_q=2$, then $q(e_{i_j}+\ldots+e_{i_l})=-1$. 
In both cases $q$ is not positive definite. 

Therefore $\langle e_{i_j},e_{i_l}\rangle_q=1$.
Note that in this case
$(i_j,i_{j+1},\ldots,i_l)$ is a~positive chordless cycle and
$[v_{i_l},v_{i_{l-1}},\ldots,v_{i_{j+1}}]$ is a~root. If $l=s+1$, then
$$v\equiv [v_{i_2},\ldots,v_{i_j},\ldots,v_{i_l}]\in\mathfrak{j}, $$
by the definition. Therefore we can assume that $l<s+1$. If $j=2$, then
$$\begin{array}{lcl}[v_{i_1},v]&\equiv& [v_{i_j},[[v_{i_l},v_{i_{l-1}},\ldots,v_{i_{j+1}}],
[v_{i_{l+1}},\ldots ,v_{i_{s+1}}]] \\
&\equiv& [[v_{i_l},v_{i_{l-1}},\ldots,v_{i_{j+1}}],[v_{i_j}
[v_{i_{l+1}},\ldots ,v_{i_{s+1}}]],
\end{array}
 $$ because $[v_{i_j},[v_{i_l},v_{i_{l-1}},\ldots,v_{i_{j+1}}]]\in\mathfrak{j}$.
It follows by 1), that $[v_{i_1},v]\in\mathfrak{j}$, because
$\langle e_{i_1},e_{i_{j+1}}+\ldots+e_{i_l}\rangle_q=0$.
Therefore we can assume that $2<j<l<s+1$ and
$$[v_{i_1},v]\equiv [v_{i_1},[x,[v_{i_j},[b,y]]]], $$ where
$x=[v_{i_{j-1}},\ldots,v_{i_2}]$, $b=[v_{i_l},v_{i_{l-1}},\ldots,v_{i_{j+1}}]$
and $y=[v_{i_{l+1}},\ldots,v_{i_{s+1}}]$.
Since $e_x$, $e_{i_j}+e_b+e_y$ and $e_x+e_{i_j}+e_b+e_y$ are roots of $q$,
by Lemma \ref{lem:pierwiastki}(b$'$) we have $\langle e_x,e_{i_j}+e_b+e_y\rangle_q=-1$.
Consider
$$
\begin{array}{lcl}
 q(-e_y+(e_{i_1}+e_x))&=&q(e_y)+q(e_{i_1}+e_x)-\langle e_{i_1}+e_x,e_y\rangle_q \\
 &=&1+1-\langle e_{i_1},e_y\rangle_q-\langle e_x,e_y\rangle_q \\
 &=&1-\langle e_x,e_y\rangle_q,
\end{array}
$$
because $\langle e_{i_1},e_y\rangle_q=\langle e_{i_1},e_{i_{s+1}}\rangle_q=1$ and
$\langle e_{i_1},e_x\rangle_q=\langle e_{i_1},e_{i_2}\rangle_q=-1$.

On the other hand
$$[v_{i_1},v]\equiv [v_{i_1},[x,[v_{i_j},[b,y]]]]=-[v_{i_1},[v_{i_j},[[b,y],x]]]-[v_{i_1},[[b,y],[x,v_{i_j}]]]. $$
Therefore $\langle e_{i_j},e_x\rangle_q=-1$, because otherwise by the induction hypothesis we have
$[x,v_{i_j}]\in\mathfrak{j}$, and by 1),
$$ [v_{i_1},v]\equiv [v_{i_1},[x,[v_{i_j},[b,y]]]]=-[v_{i_1},[v_{i_j},[[b,y],x]]]\in\mathfrak{j}.$$
Similarly we have $\langle e_x,e_b\rangle_q=-1$, because otherwise
$$\begin{array}{lcl}[v_{i_1},v]&\equiv&[v_{i_1},[x,[v_{i_j},[b,y]]]] 
\equiv[v_{i_1},[x,[b,[v_{i_j},y]]]] 
\equiv [v_{i_1},[b,[[y,v_{i_j}],x]]]\in \mathfrak{j},
  \end{array}
$$ by the case 1).

Finally $$\langle e_x,e_y\rangle_q=
-1-\langle e_x,e_{i_j}+e_b\rangle_q=-1-\langle e_x,e_{i_j}\rangle_q-\langle e_x,e_b\rangle_q=1$$
and $q(-e_y+(e_x+e_{i_1}))=1-\langle e_x,e_y\rangle_q=0$. This is a~contradiction, because
$q$ is positive definite. This finishes the proof.
\epv

\section{Examples and final remarks}\label{sec:lie-rep-dir}

In this section we present some examples and remarks that illustrate
basic results of this paper.

\begin{thm}\punkt
 If $A=\mathbb{C}Q/I$ is a~representation directed $\mathbb{C}$-algebra,
such that its Tits form $q_A$
is positive definite, then the map
\begin{equation}\Phi:L(q_A,\mathfrak{j})\to\CK(A)
 \end{equation}
 given by $v_i\mapsto u_i$ is an~isomorphism of Lie algebras.
Moreover $L(q_A,\mathfrak{j})\cong G^+(q_A)$.
\label{thm:hall-rel}\end{thm}

\textbf{Proof.} By Corollary \ref{cor:iso-lie}, the map
$\Phi:L(q_A,\mathfrak{r})\to\CK(A)$, given by $\Phi(v_i)=u_i$, is an~isomorphism of
Lie algebras. By Propositions \ref{prop:first-red} and \ref{prop:second-reduction}, we have
$L(q_A,\mathfrak{j})=L(q_A,\mathfrak{r})$, because $q_A$ is positive definite.
The isomorphism $L(q_A,\mathfrak{j})\cong G^+(q_A)$ follows from Proposition \ref{prop:nonneg}.
\epv

\begin{rem}\punkt {\rm Let $A$ be a~representation directed $\mathbb{C}$-algebra and let $q_A$ be its Tits form.
It is well-known (see \cite{bo83}) that $q_A$ is weakly positive.
It follows that the set $\CR_{q_A}^+$ of positive roots of $q_A$
is finite. Therefore $\dim_{\mathbb{C}}\CK(A)=|\CR_{q_A}^+|$ is
finite. In this case, the subset $\mathfrak{r}$ of $L(q_A)$ is
finite, even if $q_A$ is not positive definite. Moreover, we are able to describe
an~algorithm that constructs the set $\mathfrak{r}$. Indeed, it is enough
to develop Definition \ref{def:lie} and construct all Weyl roots of $q$ (see \cite[Remark 4.15]{kos01}).

If $q_A$ is positive definite, then $L(q_A,\mathfrak{j})=L(q_A,\mathfrak{r})$. The set $\mathfrak{j}$
is a~minimal set generating the ideal $(\mathfrak{r})$ and $\mathfrak{j}$
is smaller than $\mathfrak{r}$ (see Example \ref{ex:relations}).

If $q_A$ is not positive definite, then $(\mathfrak{j})\subsetneq (\mathfrak{r})$
in general (see Example \ref{ex:non-positive}).}
\end{rem}

\begin{ex}\punkt {\rm Let $L$ be the~following poset
$$\xymatrix{&&4&\\
L:&2\ar[ur]&&3\ar[ul]\\
&&1\ar[ur]\ar[ul]&}$$ and let $KL$ be the incidence algebra of the
poset $L$ (see \cite{Si92}). It is easy to see that $KL$ is representation directed,
$q_{KL}$ is positive definite and $B(q_{KL})$ has
the form
$$
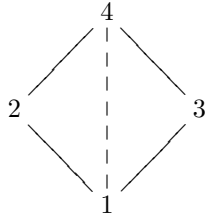
\xymatrix{&&4&\\
&2\ar@{-}[ur]&&3\ar@{-}[ul]\\
&&1\ar@{-}[ur]\ar@{-}[ul]\ar@{--}[uu]&}$$ Then
$$\begin{array}{lcl}\mathfrak{j}&=&\{[u_2,u_3],[u_1,u_4],[u_1,[u_1,u_2]],[u_1,[u_1,u_3]],
[u_2,[u_1,u_2]],[u_3,[u_1,u_3]],\\ && [u_2,[u_2,u_4]],
[u_3,[u_3,u_4]], [u_4,[u_2,u_4]],[u_4,[u_3,u_4]],[u_1[u_2,u_4]],
\\ &&[u_1,[u_3,u_4]]\}
\end{array}
$$
and $L(q_{KL},\mathfrak{j})\cong\CK(KL)$.
Note that \small
 $$\begin{array}{lcl}\mathfrak{r}&=&\{[u_2,u_3],[u_1,u_4],[u_1,[u_1,u_2]],[u_1,[u_1,u_3]],
[u_2,[u_2,u_1]],[u_3,[u_3,u_1]],\\&&  [u_2,[u_2,u_4]],  [u_3,[u_3,u_4]],
[u_4,[u_4,u_2]],[u_4,[u_4,u_3]],[u_1[u_2,u_4]],\\ &&[u_1,[u_3,u_4]],
[u_3,u_2],[u_4,u_1],[u_1,[u_2,u_1]],[u_1,[u_3,u_1]],
[u_2,[u_1,u_2]],\\ &&[u_3,[u_1,u_3]],[u_2,[u_4,u_2]],  [u_3,[u_4,u_3]],
[u_4,[u_2,u_4]],[u_4,[u_3,u_4]],\\ &&[u_1[u_4,u_2]],[u_1,[u_4,u_3]],
[u_1,u_2,u_1,u_3], [u_2,u_2,u_1,u_3], [u_3,u_2,u_1,u_3], \\ &&[u_1,u_3,u_1,u_2], [u_2,u_3,u_1,u_2], [u_3,u_3,u_1,u_2],
[u_1,u_2,u_3,u_1],\\ && [u_2,u_2,u_3,u_1], [u_3,u_2,u_3,u_1], [u_1,u_3,u_2,u_1], [u_2,u_3,u_2,u_1], \\ &&[u_3,u_3,u_2,u_1],
[u_4,u_2,u_4,u_3], [u_2,u_2,u_4,u_3], [u_3,u_2,u_4,u_3],\\ && [u_4,u_3,u_4,u_2], [u_2,u_3,u_4,u_2], [u_3,u_3,u_4,u_2],
[u_4,u_2,u_3,u_4],\\ && [u_2,u_2,u_3,u_4], [u_3,u_2,u_3,u_4], [u_4,u_3,u_2,u_4], [u_2,u_3,u_2,u_4], [u_3,u_3,u_2,u_4],
\\ && [u_1,u_4,u_2,u_3,u_1],
[u_2,u_4,u_2,u_3,u_1], [u_3,u_4,u_2,u_3,u_1], [u_4,u_4,u_2,u_3,u_1],
\\ && [u_1,u_4,u_3,u_2,u_1],
[u_2,u_4,u_3,u_2,u_1], [u_3,u_4,u_3,u_2,u_1], [u_4,u_4,u_3,u_2,u_1],
\\ && [u_1,u_1,u_2,u_3,u_4],
[u_2,u_1,u_2,u_3,u_4], [u_3,u_1,u_2,u_3,u_4], [u_4,u_1,u_2,u_3,u_4],
\\ && [u_1,u_1,u_3,u_2,u_4],
[u_2,u_1,u_3,u_2,u_4], [u_3,u_1,u_3,u_2,u_4], [u_4,u_1,u_3,u_2,u_4],
\\ && [u_1,u_4,u_2,u_1,u_3],
[u_2,u_4,u_2,u_1,u_3], [u_3,u_4,u_2,u_1,u_3], [u_4,u_4,u_2,u_1,u_3],
\\ && [u_1,u_4,u_3,u_1,u_2],
[u_2,u_4,u_3,u_1,u_2], [u_3,u_4,u_3,u_1,u_2], [u_4,u_4,u_3,u_1,u_2],
\\ && [u_1,u_1,u_2,u_4,u_3],
[u_2,u_1,u_2,u_4,u_3], [u_3,u_1,u_2,u_4,u_3], [u_4,u_1,u_2,u_4,u_3],
\\ && [u_1,u_1,u_3,u_4,u_2],
[u_2,u_1,u_3,u_4,u_2], [u_3,u_1,u_3,u_4,u_2], [u_4,u_1,u_3,u_4,u_2]\}
\end{array}
$$\normalsize } \label{ex:relations}\end{ex}

\begin{ex}\punkt
 {\rm Consider the following graph $$
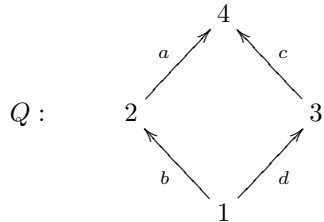
\xymatrix{&&4&\\
Q:&2\ar[ur]^a&&3\ar[ul]_c\\
&&1\ar[ur]_d\ar[ul]^b&}$$
Let $A=\mathbb{C}Q/I$, where $I=(ab,cd)$. The form $q_A$ is not positive definite
and $B(q_A)$ has the following form $$
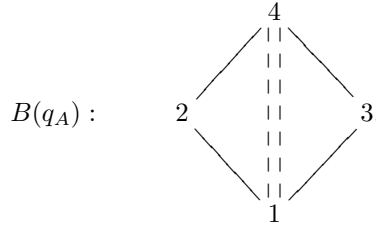
\xymatrix{&&4&\\
B(q_A):&2\ar@{-}[ur]&&3.\ar@{-}[ul]\\
&&1\ar@{-}[ur]\ar@{-}[ul]\ar@{--}@<0.5ex>[uu]\ar@{--}@<-0.5ex>[uu]&}$$
Note that $[u_4,[u_3,[u_2,u_1]]]\in(\mathfrak{r})$, but
$[u_4,[u_3,[u_2,u_1]]]\not\in(\mathfrak{j})$. 
On the other hand, the algebra $A$ is representation directed and $q_A$
is weakly positive. By Corollary \ref{cor:iso-lie}, $\CK(A)\cong L(q_A,\mathfrak{r})$.}
\label{ex:non-positive}\end{ex}

\end{document}